\documentclass[11pt]{article}

\usepackage{amsfonts}
\usepackage{amssymb,amsmath,amsthm}
\usepackage{latexsym}
\usepackage{fullpage}
\usepackage{hyperref}
\usepackage{graphicx}
\usepackage{tkz-euclide,subfigure}
\usepackage{tikz}
\usepackage[utf8]{inputenc}
\usepackage{pgfplots}
\usetikzlibrary{shapes.misc}
\usepackage{ amssymb }

\newtheorem{theorem}{Theorem}[section]

\newtheorem{proposition}[theorem]{Proposition}

\newtheorem{lemma}[theorem]{Lemma}
\newtheorem{cor}[theorem]{Corollary}
\newtheorem{definition}[theorem]{Definition}

\newtheorem{example}[theorem]{Example}

\theoremstyle{definition}
\newtheorem{remark}[theorem]{Remark}

\newcounter{tenumerate}

\renewcommand{\epsilon}{\varepsilon}

\newcommand{\remove}[1]{}

\renewcommand{\leq}{\leqslant}
\renewcommand{\geq}{\geqslant}

\def\XXint#1#2#3{{\setbox0=\hbox{$#1{#2#3}{\int}$}
\vcenter{\hbox{$#2#3$}}\kern-.5\wd0}}

\title{\textbf{Well-posedness of quadratic RBSDEs and BSDEs with one-sided growth restrictions}}
\begin{document}
\author{{Shiqiu Zheng\thanks{E-mail: shiqiu@tust.edu.cn}}\\
  \\
\small College of Sciences, Tianjin University of Science and Technology, Tianjin 300457, China}
\date{}
\maketitle
\begin{abstract}
In this paper, we investigate the well-posedness of bounded and unbounded solutions for reflected backward stochastic differential equations (RBSDEs) and backward stochastic differential equations (BSDEs). The generators of these equations satisfy a one-sided growth restriction on the variable $y$ and have a general quadratic growth in the variable $z$. The solutions $Y_t$ (and the obstacles for RBSDEs) take values in either $\mathbf{R}$ or $(0, \infty)$. We obtain the existence of solutions primarily by using the methods from Essaky and Hassani (2011) and Bahlali et al. (2017). For the uniqueness of solutions, we provide a method applicable when the generators are convex in $(y,z)$ or are (locally) Lipschitz in $y$ and convex in $z$. Our method relies on the $\theta$-difference technique introduced by Briand and Hu (2008), and some novel comparison arguments based on RBSDEs. We also establish some general comparison theorems for such RBSDEs and BSDEs.\\\\
\textbf{Keywords:} reflected backward stochastic differential equation; backward stochastic differential equation; comparison theorem; quadratic growth; one-sided growth\\
\textbf{AMS Subject Classification:} 60H10
\end{abstract}
\section{Introduction}
A backward stochastic differential equation (BSDE) is usually called a quadratic BSDE, when its generator $g(t,y,z)$ has a quadratic growth in the variable $z$. There have been many studies on the well-posedness of one-dimensional quadratic BSDEs. We refer to \cite{K, LS, BL, BE, F, Tian, Im} for the well-posedness of bounded solutions, and to \cite{BH06, BH08, D, DR, EH, BK, F, B17, B19, BT, Y, FHT1, FHT2, LG} for the well-posedness of unbounded solutions. Furthermore, for more results on one-dimensional BSDEs and their applications, we refer to a user's guide recently provided by Fan et al. \cite{FHT2}. A reflected BSDE (RBSDE) can be considered as a BSDE containing an additional nondecreasing process $K$ to push the solution $Y$ above a given obstacle $L$ under the Skorokhod condition. The theory of RBSDEs was initiated by El Karoui et al. \cite{EPP} and has important applications in several areas, including the pricing of American options, mixed stochastic control, and obstacle problem for PDEs (see \cite{EPP, K2, Xu, Li, BY, EH, Po, Gu} and the references therein). The extension of well-posedness results from quadratic BSDEs to the RBSDE setting has also been investigated. We refer to \cite{K2, Xu, Li} for the well-posedness of bounded solutions, and to \cite{BY, EH, Gu} for the well-posedness of unbounded solutions.

In this paper, we investigate the well-posedness of quadratic RBSDEs and BSDEs, whose generators $g(t,y,z)$ have a general growth in $y$ and a general quadratic growth in $z$ (see Assumption (\hyperref[2A1]{2A1})), and solutions $Y$ (and obstacles $L$ for RBSDEs) take values in $D=\bf{R}$ or $(0,\infty)$. We first consider the comparison theorems for such quadratic RBSDEs and BSDEs. We show that for two RBSDEs, the comparison of the solutions $Y$ implies the comparison of the solutions $K$ (see Proposition \ref{pr3.1}). This result was then used to obtain an existence result of maximal solutions of such RBSDEs (resp. BSDEs) (see Proposition \ref{pr3.3}), which plays a crucial role in the study of the uniqueness of solutions in this paper. We also establish some general comparison theorems for the solutions $Y$, which roughly indicate that the solutions $Y$ can be compared whenever the solutions of such RBSDEs (resp. BSDEs) are unique in a certain space (see Remark \ref{r3.5}). This phenomenon was previously observed by Zheng \cite[Theorem 2.7]{Zheng1} for BSDEs whose
generators are Lipschitz in $y$, in a general setting.

To get the well-posedness of bounded solutions of an RBSDE under (\hyperref[2A1]{2A1}), we further assume that the generator satisfies
\begin{equation}\label{1.1}
1_{\{y\geq c\}}g(t,y,0)\leq u(t)l(y),\quad\forall y\in D,\tag{1.1}
\end{equation}
for a constant $c>0$, a nonnegative integrable function $u(t)$ on $[0,T]$ and a positive continuous function $l(y)$ satisfying $\int_{-\infty}^01/l(y)dy=\int_0^\infty1/l(y)dy=\infty$. To get the well-posedness of unbounded solutions of an RBSDE under (\hyperref[2A1]{2A1}), we further assume that the generator satisfies
\begin{equation}\label{1.2}
1_{\{y\geq c\}}g(t,y,z)\leq \delta_t +\gamma_t |y|+\kappa|z|+f(y)|z|^2, \quad\forall (y,z)\in D\times\mathbf{R}^d, \tag{1.2}
\end{equation}
for two constants $c>0, \kappa\geq0$, two nonnegative processes $\delta,\gamma$, and a nonnegative continuous function $f(y)$ on $D$. To get the well-posedness of bounded (resp. unbounded) solutions of a BSDE under (\hyperref[2A1]{2A1}), when $D=\mathbf{R}$, we further assume that the generator satisfies (\ref{1.1}) (resp. (\ref{1.2})) and
\begin{equation}\label{1.3}
1_{\{y\leq -c\}}g(t,y,0)\geq-u(t)l(y),\quad \forall y\in D,\tag{1.3}
\end{equation}
\begin{equation}\label{1.4}
(\textrm{resp.} \quad 1_{\{y\leq -c\}}g(t,y,z)\geq -\delta_t -\gamma_t |y|-\kappa|z|-f(-y)|z|^2,\quad\forall(y,z)\in D\times\mathbf{R}^d), \tag{1.4}
\end{equation}
and when $D=(0,\infty)$, we further assume that the generator satisfies (\ref{1.1}) (resp. (\ref{1.2})) and
\begin{equation}\label{1.5}
1_{\{y\leq b\}}g(t,y,0)\geq -u(t)|y\ln(y)|,\quad \forall y\in D,\tag{1.5}
\end{equation}
\begin{equation}\label{1.6}
(\textrm{resp.}\quad 1_{\{y\leq b\}}g(t,y,z)\geq -\vartheta_t|y|-\kappa|z|-\frac{\nu}{y}|z|^2, \quad \forall (y,z)\in D\times\mathbf{R}^d),\tag{1.6}
\end{equation}
for two constants $b>0$, $\nu\geq0$ and a nonnegative process $\vartheta$.

A key difference between the quadratic growth conditions in (\ref{1.1})-(\ref{1.6}) and those in existing studies lies in the one-sided growth in $y$. These one-sided growth conditions contain singular generators and generators with general stochastic coefficients (see Examples \ref{ex4.6} and \ref{ex5.10}). Specifically:

\begin{itemize}
  \item (\ref{1.1}) and (\ref{1.3}) indicate that the generator has a one-sided superlinear growth in $y$ and a general quadratic growth in $z$. This extends the quadratic growth conditions for bounded solutions in \cite{K, LS, K2, BL, Xu, BE, F, Tian, Im}. Furthermore, when $u(t)$ is continuous and the terminal time $T$ is small enough, $l(y)$ in (\ref{1.1}) and (\ref{1.3}) can be an arbitrary nonnegative continuous function and $D$ can be an arbitrary open interval. This includes the ``characteristic BSDE" observed by Ma et al. \cite[Equation (3.8)]{MZ} in the study of the well-posedness of forward BSDEs. (\ref{1.2}) and (\ref{1.4}) are more general than the one-sided conditions for quadratic BSDEs assumed in \cite{FHT1, FHT2} (see Remark \ref{r5.1}(ii)).
  \item When $D=(0,\infty)$, the one-sided conditions imply that the generator $g(t,\cdot,z)$ can be singular at $0$. Recently, such singular BSDEs and RBSDEs have been studied in several cases. For example, the BSDEs with nonnegative generators whose quadratic growth term takes the form $|z|^2/y$ were studied by \cite{BT, LG} (see Remark \ref{r5.8} for a contrast); BSDEs and RBSDEs with some special generators containing the term $f(y)|z|^2$ were studied by \cite{Zheng2, Zheng3}, where $f$ is a locally integrable function defined on an open interval.
\end{itemize}
Some special cases of the BSDEs and RBSDEs discussed in this paper have important applications, for instance, in economics and finance (see \cite{DE, SS, BT, Tian, LG, Zheng3}), in the probabilistic interpretation of viscosity solutions to PDEs with singular coefficients (see \cite{Zheng2, Zheng3}), and in the study of the well-posedness of forward BSDEs (see \cite[Equation (3.8) and Remark 3.1]{MZ}). These applications are the main motivation for our current investigation into the general case of these equations (see Example \ref{ex4.5} for a motivating example from finance).

We establish the existence of solutions primarily using the methods inspired by Essaky and Hassani \cite{EH} and Bahlali et al. \cite{B17}, as well as some well-posedness results for ODEs. However, it seems to be difficult to prove the uniqueness of solutions due to the singularity and general growth properties of the generators. To address this, this paper provides a method applicable when their generators satisfy a $\theta$-domination condition (see (\ref{2.1})). This $\theta$-domination condition is inspired by Fan and Hu \cite[Assumption (H2')]{FH} (see also \cite{FHT1, FHT2}). It includes generators that are convex in $(y,z)$, generators that are (locally) Lipschitz in $y$ and convex in $z$, as well as some non-convex generators. Our method relies on the $\theta$-difference technique introduced by Briand and Hu \cite{BH08}, and some novel comparison arguments based on RBSDEs. It differs from the methods used for quadratic RBSDEs in \cite{K2, BY, Li, Gu}. The use of the comparison arguments is a key difference between our method and those based on the $\theta$-difference technique for quadratic BSDEs in \cite{BH08, Y, FH, FHT1, FHT2, LG}. For the bounded solutions of RBSDEs, our method can be described briefly as follows:
\begin{itemize}
\item We first prove that the RBSDE$(g,\xi,L)$ under (\ref{1.1}) admits a minimal solution $(Y,Z,K)$ such that the range of $Y$ is included in a closed subset of $D$.
\item Then, we prove that under (\ref{2.1}) and (\hyperref[4A2]{4A2}), for any solution $(\hat{Y},\hat{Z},\hat{K})$ to the RBSDE$(g,\xi,L)$ such that the range of $\hat{Y}$ is included in a closed subset of $D$, by using some comparison arguments based on RBSDEs, we can find an RBSDE, which admits a maximal solution $(\bar{Y},\bar{Z},\bar{K})$ such that for each $\theta\in(0,1)$, $\bar{Y}_t\geq\frac{\hat{Y}_t-\theta Y_t}{1-\theta}.$
\item Finally, since $\bar{Y}_t\geq\frac{\hat{Y}_t-\theta Y_t}{1-\theta}$, when $\theta$ tends to $1$, we get $\hat{Y}_t\leq Y_t$, i.e., $(Y_t,Z_t,K_t)=(\hat{Y}_t,\hat{Z}_t,\hat{K}_t).$
\end{itemize}
The BSDEs case follows the spirit of the method above, but requires some different treatments. It is known that a locally Lipschitz condition (which does not necessarily satisfy the $\theta$-domination condition: (\ref{2.1}) together with (\hyperref[4A2]{4A2})) typically guarantees the uniqueness of bounded solutions (see \cite{K2, BE, Li, Im}). Since the generators in our framework may fail to satisfy such a condition (see Example \ref{ex4.6}(ii)), to ensure the uniqueness of bounded solutions, we also impose the $\theta$-domination condition.

This paper is organized as follows. In Section 2, we present some assumptions and lemmas. In Section 3, we study comparison theorems. In Sections 4 and 5, we study the well-posedness of bounded solutions and unbounded solutions, respectively. In the Appendix, we present some auxiliary results.
\section{Preliminaries}
Let $(\Omega ,\mathcal{F},\mathit{P})$ be a complete probability space, on which a $d$-dimensional standard Brownian motion $(B_t)_{t\geq
0}$ is defined. Let $({\mathcal{F}}_t)_{t\geq 0}$ be
the natural filtration generated by $(B_t)_{t\geq 0}$, augmented
by the $\mathit{P}$-null sets of ${\mathcal{F}}$. Let $|z|$ denote the
Euclidean norm of $z\in {\mathbf{R}}^d$. Let $x\cdot z$ denote the scalar product of $x,z\in\mathbf{R}^d$. Let ${\mathcal{P}}$ be the progressive measurable sigma-field on $[0,T]\times\Omega.$ Let $T>0$ and $p>1$ be given real numbers. We always assume that $D=\bf{R}$ or $D=(0,\infty)$, except for Subsection \ref{sec4.2}, where $D$ is an arbitrary open interval. We introduce the following spaces:

$C_+(D):=\{f: D\rightarrow \mathbf{R}$, nonnegative and continuous$\}$;

$L_{loc}(D):=\{f: D\rightarrow \mathbf{R},$ measurable and locally integrable$\}$;

$L^1_+[0,T]:=\{f:[0,T]\rightarrow \mathbf{R},$ measurable, nonnegative and integrable$\}$;

$L_D({\mathcal{F}}_T):=\{\xi: {\mathcal{F}}_T$-measurable random variable taking values in $D\}$;

$L^r_D({\mathcal {F}}_T):=\{\xi\in L_D({\mathcal{F}}_T):$ ${{E}}\left[|\xi|^r\right]<\infty\},\ r\geq1; $

$L^\infty_D({\mathcal {F}}_T):=\{\xi\in L({\mathcal{F}}_T):$ $\xi$ takes values in a bounded closed subset of $D\};$

${\mathcal{C}}_D:=\{(\psi_t)_{t\in[0,T]}:$ continuous and $({\mathcal{F}}_t)$-adapted process taking values in $D\}$;

${\mathcal{S}}^r_D:=\{(\psi_t)_{t\in[0,T]}:$ process in ${\mathcal{C}}_D$ such that ${{E}} \left
[{\mathrm{sup}}_{0\leq t\leq T} |\psi _t|^r\right]<\infty \},\ r\geq1;$

${\mathcal{S}}_D^\infty:=\{(\psi_t)_{t\in[0,T]}:$ process in ${\mathcal{C}}$ taking values in a bounded closed subset of $D\};$

${\mathcal{A}}:=\{(\psi_t)_{t\in[0,T]}:$ increasing, continuous, $({\mathcal{F}}_t)$-adapted
$\mathbf{R}$-valued process with $\psi_0=0$$\}$;

$H^r_d:=\{(\psi_t)_{t\in[0,T]}:$  $\mathbf{R}^d$-valued, $({\mathcal{F}}_t)$-progressively measurable and $\int_0^T|\psi_t|^rdt
<\infty\},\ r\geq1;$

${\mathcal{H}}_d^r:=\{(\psi_t)_{t\in[0,T]}:$ process in $H^2_d$ such that ${{E}}[(\int_0^T|\psi_t|^2dt)^{\frac{r}{2}}]<\infty \},\ r\geq1;$

${\mathcal{H}}_d^{BMO}:=\{(\psi_t)_{t\in[0,T]}:$ process in $H_d^2$ such that $\sup_{\tau\in{\mathcal{T}}_{0,T}}\|E[\int_\tau^T|\psi_t|^2dt|{\mathcal{F}}_\tau]\|_\infty<\infty\}.$\\
For convenience, when $D=\mathbf{R}$ is clear, we write $L_D(\mathcal{F}_T)$, $\mathcal{C}_D$ and $\mathcal{S}_D^r$ as $L(\mathcal{F}_T)$, $\mathcal{C}$ and $\mathcal{S}^r$, respectively. Note that in this paper, all the equalities and inequalities for random variables are understood in the almost sure sense.

Throughout, we assume that $g$ is a function which satisfies the following assumption:
\begin{itemize}
\item \textbf{(2A1)}\label{2A1}
\begin{equation*}
g\left(\omega,t,y,z\right): \Omega \times [0,T]\times D\times\mathbf{R}^d\longmapsto \mathbf{R},
\end{equation*}
is measurable with respect to ${\mathcal{P}}\otimes({\mathcal{B}}({\mathbf{R}})\cap D)\otimes{\mathcal{B}}({\mathbf{R}}^d)$ and satisfies the following two conditions:
\begin{itemize}
\item[(i)] $dt\times dP\textrm{-}a.e.,$ $g(t,\cdot,\cdot)$ is continuous on $D\times\mathbf{R}^{\mathit{d}}$;
\item[(ii)] for any two processes $X,Y\in{\mathcal{C}}_{D}$ satisfying $X_t\leq Y_t$ for all $t\in[0,T]$, there exist two nonnegative processes $C\in{\mathcal{C}}$ and $\eta\in H^1_1$ such that $dt\times dP\textrm{-}a.e.,$ for each $y\in[X_t(\omega),Y_t(\omega)]$ and $z\in\mathbf{R}^d$, $$|g(t,y,z)|\leq\eta_t+C_t|z|^2.$$
\end{itemize}
\end{itemize}

\begin{remark}\label{r2.1}
(\hyperref[2A1]{2A1})(ii) implies that $g(t,y,z)$ has a general growth in $y$ and a general quadratic growth in $z$. In fact, by Essaky and
Hassani \cite[Remark 2.2(3)]{EH1}, we get that (\hyperref[2A1]{2A1})(ii) is satisfied, if there exist a nonnegative process $\eta\in H^1_1$ and two functions $\phi(\omega,t,y),\psi(\omega,t,y): \Omega\times [0,T]\times D\longmapsto[0,\infty)$ which are both measurable with respect to ${\mathcal{P}}\otimes({\mathcal{B}}({\mathbf{R}})\cap D)$ and continuous on $[0,T]\times D$, such that $dt\times dP\textrm{-}a.e.,$ for each $(y,z)\in D\times{\mathbf{ R}}^{\mathit{d}},$
\begin{equation*}
|g(t,y,z)|\leq \eta_t+\phi(t,y)+\psi(t,y)|z|^2.
\end{equation*}
\end{remark}

To study the uniqueness of solutions, we introduce the following $\theta$-domination condition, which is inspired by \cite{FH}, and includes generators that are convex in $(y,z)$, generators that are (locally) Lipschitz in $y$ and convex in $z$, as well as some non-convex generators (see Remark \ref{r4.1}).

\begin{itemize}
  \item\textbf{($\theta$-domination condition)} We say that $g$ satisfies the \textbf{$\theta$-domination condition for ${\mathcal{C}}_D$ and $H^2_d$ with $G$}, if there exists a function $G(\omega,t,y,z,x_1,x_2): \Omega \times [0,T]\times D\times\mathbf{R}^{\mathit{d}}\times D\times D\longmapsto\mathbf{R},$ which is measurable with respect to ${\mathcal{P}}\otimes({\mathcal{B}}({\mathbf{R}})\cap D)\otimes{\mathcal{B}}({\mathbf{R}}^d)\otimes({\mathcal{B}}({\mathbf{R}})\cap D)\otimes({\mathcal{B}}({\mathbf{R}})\cap D)$ such that for any $y^1,y^2\in{\mathcal{C}}_D$ and $z^1,z^2\in H^2_d$, $G(\omega,t,y,z,y^1_t,y^2_t)$ satisfies (\hyperref[2A1]{2A1}) and $dt\times dP\textrm{-}a.e.,$ for each $\theta\in(0,1)$, if the range of $\frac{y_t^1-\theta y_t^2}{1-\theta}$ is included in $D$, then
\begin{equation*}\label{2.1}
{g}(t,y^1_t,z_t^1)-\theta{g}(t,y^2_t,z^2_t)\leq (1-\theta)G\left(t,\frac{y_t^1-\theta y_t^2}{1-\theta},\frac{z_t^1-\theta z_t^2}{1-\theta},y^1_t,y^2_t\right).\tag{2.1}
\end{equation*}
\end{itemize}

Let $\xi\in L_D({\mathcal{F}}_T)$ and $L\in{\mathcal{C}}_D$ be given such that $\xi\geq L_T$. We consider the following RBSDE:
\begin{equation*}\label{2.2}
\left\{
   \begin{array}{ll}
    Y_t=\xi +\int_t^Tg(s,Y_s,Z_s)ds+K_T-K_t-\int_t^TZ_s\cdot dB_s,\ \ t\in[0,T],\\
    \forall t\in[0,T],\ \ Y_t\geq L_t,\\
    \int_0^T(Y_t-L_t)dK_t=0,\ \ \ (\textrm{Skorokhod condition})
   \end{array}
 \right.
\tag{2.2}
\end{equation*}
and BSDE:
\begin{equation*}\label{2.3}
Y_t=\xi +\int_t^Tg(s,Y_s,Z_s)ds-\int_t^TZ_s\cdot dB_s,\ \ t\in[0,T],\tag{2.3}
\end{equation*}
where $T$ is the terminal time, $\xi$ is the terminal variable, $g$ is the generator, and $L$ is the lower obstacle. We denote (\ref{2.2}) and (\ref{2.3}) by the RBSDE$(g,\xi,L)$ and the BSDE$(g,\xi)$, respectively.

\begin{definition} A solution of the RBSDE$(g,\xi,L)$ is a triple of processes $(Y,Z,K)\in{\mathcal{C}}_D
\times{H}^2_d\times{\mathcal{A}},$ which satisfies $\int_0^{T}|g(s,Y_s,Z_s)|ds<\infty$ and (\ref{2.2}). A solution of the BSDE$(g,\xi)$ is a pair of processes $(Y,Z)\in{\mathcal{C}}_D
\times{H}^2_d,$ which satisfies $\int_0^{T}|g(s,Y_s,Z_s)|ds<\infty$ and (\ref{2.3}).
\end{definition}

Note that we say that $(Y,Z,K)$ is a \textbf{unique (resp. minimal or maximal) solution} of the RBSDE$(g,\xi,L)$ such that $Y$ satisfies some condition (\textbf{C}), if it is a solution of the RBSDE$(g,\xi,L)$ such that $Y$ satisfies the condition (\textbf{C}), and for each solution $(Y',Z',K')$ of the RBSDE$(g,\xi,L)$ such that $Y'$ satisfies the condition (\textbf{C}), we have $Y_t=Y'_t$ (resp. $Y_t\leq Y'_t$ or $Y_t\geq Y'_t$) for each $t\in[0,T].$ The BSDEs case is similar.

In this section, we introduce two semimartingales:
\begin{equation*}\label{2.4}
Y^1_t=Y^1_T+\int_t^Th_1(s)ds+A^1_T-A^1_t-\int_t^TZ^1_s\cdot dB_s,\ \ t\in[0,T],\tag{2.4}
\end{equation*}
and
\begin{equation*}\label{2.5}
Y^2_t=Y^2_T+\int_t^Th_2(s)ds-A^2_T+A^2_t-\int_t^TZ^2_s\cdot dB_s,\ \ t\in[0,T],\tag{2.5}
\end{equation*}
where $Y^i\in{\mathcal{C}}_D$, $A^i\in {\mathcal{A}},$ $Z^i\in H^2_d$ and $h_i\in H^1_1, i=1,2.$  We say that \textbf{the RBSDE$(g,\xi,L)$ is dominated by $Y^1$}, if the following (i) and (ii) hold:

(i) $\xi\leq Y^1_T$ and for all $t\in[0,T]$, $L_t\leq Y^1_t$;

(ii) $g(t,Y_t^1,Z_t^1)\leq h_1(t),\ dt\times dP\textrm{-}a.e.$\\
We say that \textbf{the BSDE$(g,\xi)$ is dominated by $Y^2$ and $Y^1$}, if the following (i) and (ii) hold:

(i) $Y^2_T\leq\xi\leq Y^1_T$ and for all $t\in[0,T]$, $Y^2_t\leq Y^1_t$;

(ii) $g(t,Y_t^2,Z_t^2)\geq h_2(t)$ and $g(t,Y_t^1,Z_t^1)\leq h_1(t),\ dt\times dP\textrm{-}a.e.$\\

In the following two lemmas, the existence of solutions follows from \cite[Theorem 3.1]{EH} or from the domination arguments in \cite{B17, B19}. The minimality of the solution in Lemma \ref{l2.3} was previously pointed out in \cite[Remark 4.2]{EH}. The proofs of Lemmas \ref{l2.3} and \ref{l2.4} are provided in Appendix \ref{appA}.
\begin{lemma}\label{l2.3}
Let the RBSDE$(g,\xi,L)$ be dominated by $Y^1$. Then it admits:

(i) A minimal solution $(y,z,k)$ such that $y\in\mathcal{C}_D$;

(ii) A maximal solution $(Y,Z,K)$ such that for all $t\in[0,T]$, $Y_t\leq Y^1_t$.
\end{lemma}
\begin{lemma}\label{l2.4}
Let the BSDE$(g,\xi)$ be dominated by $Y^2$ and $Y^1$. Then it admits:

(i) A minimal solution $(y,z)$ such that for all $t\in[0,T]$, $y_t\geq Y^2_t$;

(ii) A maximal solution $(Y,Z)$ such that for all $t\in[0,T]$, $Y_t\leq Y^1_t$.\\
Moreover, for all $t\in[0,T]$, $Y_t^2\leq y_t\leq Y_t\leq Y^1_t$.
\end{lemma}
\section{Comparison theorems}\label{sec3}
We first provide a slight generalization of the comparison theorem in Lionnet \cite[Proposition 5]{Li}, using a similar method. It shows that for RBSDEs, the comparison of solutions $Y$ implies the comparison of solutions $K$.

\begin{proposition}\label{pr3.1} Let $\bar{g}$ satisfy (\hyperref[2A1]{2A1}). Let $\bar{\xi}\in L_D({\mathcal{F}}_T)$ and $\bar{L}\in{\mathcal{C}}_D$ such that $\bar{\xi}\geq\bar{L}_T$. Let the RBSDE$(\bar{g},\bar{\xi}, \bar{L})$ admit a solution $(\bar{Y},\bar{Z},\bar{K})$, and let the RBSDE$(g,\xi, L)$ admit a solution $(Y,Z,K)$. If $\bar{Y}_t\geq Y_t$ and $\bar{g}(t,\bar{Y}_t,\bar{Z}_t)\geq g(t,\bar{Y}_t,\bar{Z}_t)$, $dt\times dP\textrm{-}a.e.,$ then for each $0\leq r<t\leq T$,
\begin{equation*}
\int_r^t1_{\{Y_s\geq\bar{L}_s\}}d\bar{K}_s\leq\int_r^t1_{\{Y_s\geq\bar{L}_s\}}dK_s.
\end{equation*}
\end{proposition}

\begin{proof} Since for all $t\in[0,T]$, $\bar{Y}_t\geq Y_t$, by considering $(\bar{Y}_t-Y_t)-(\bar{Y}_t-Y_t)^+,$ we deduce that for any $0\leq r<t\leq T,$
\begin{align*}\label{3.1}
\int_r^t1_{\{\bar{Y}_s=Y_s\}}dK_s-\int_r^t1_{\{\bar{Y}_s=Y_s\}}d\bar{K}_s
=\int_r^t&1_{\{\bar{Y}_s=Y_s\}}(\bar{g}(s,\bar{Y}_s,\bar{Z}_s)-g(s,Y_s,Z_s))ds
\\&-\int_r^t1_{\{\bar{Y}_s=Y_s\}}(\bar{Z}_s-Z_s)\cdot dB_s+\frac{1}{2}\ell_t^0(\bar{Y}-Y)
-\frac{1}{2}\ell_r^0(\bar{Y}-Y),\tag{3.1}
\end{align*}
where $\ell_t^0(\bar{Y}-Y)$ is the local time of $\bar{Y}-Y$ at time $t$ and level $0$. This implies that $\int_0^T1_{\{\bar{Y}_s=Y_s\}}(\bar{Z}_s-Z_s)\cdot dB_s=0$, which leads to
\begin{equation*}\label{3.2}
1_{\{\bar{Y}_t=Y_t\}}|\bar{Z}_t-Z_t|=0,\quad dt\times dP\textrm{-}a.e.\tag{3.2}
\end{equation*}
By (\ref{3.1}) and (\ref{3.2}), we have
\begin{equation*}
\int_r^t1_{\{\bar{Y}_s=Y_s\}}dK_s-\int_r^t1_{\{\bar{Y}_s=Y_s\}}d\bar{K}_s
\geq\int_r^t1_{\{\bar{Y}_s=Y_s\}}(\bar{g}(s,\bar{Y}_s,\bar{Z}_s)-g(s,\bar{Y}_s,\bar{Z}_s))ds.
\end{equation*}
Since $\bar{g}(s,\bar{Y}_s,\bar{Z}_s)\geq g(t,\bar{Y}_t,\bar{Z}_t)$, $dt\times dP\textrm{-}a.e.,$ we have
\begin{equation*}\label{3.3}
\int_r^t1_{\{\bar{Y}_s=Y_s\}}d\bar{K}_s\leq\int_r^t1_{\{\bar{Y}_s=Y_s\}}dK_s.\tag{3.3}
\end{equation*}
Moreover, we also have
\begin{equation*}\label{3.4}
\int_0^T 1_{\{\bar{Y}_s>Y_s\}}1_{\{Y_s\geq\bar{L}_s\}}d\bar{K}_s
\leq\int_0^T1_{\{\bar{Y}_s>\bar{L}_s\}}d\bar{K}_s=0.\tag{3.4}
\end{equation*}
Then, by (\ref{3.3}) and (\ref{3.4}), we get that for any $0\leq r\leq t\leq T,$
\begin{eqnarray*}
\int_r^t1_{\{Y_s\geq\bar{L}_s\}}d\bar{K}_s&=&\int_r^t 1_{\{\bar{Y}_s>Y_s\}}1_{\{Y_s\geq\bar{L}_s\}}d\bar{K}_s
+\int_r^t1_{\{\bar{Y}_s=Y_s\}}1_{\{Y_s\geq\bar{L}_s\}}d\bar{K}_s\\
&\leq&\int_r^t1_{\{\bar{Y}_s=Y_s\}}1_{\{Y_s\geq\bar{L}_s\}}dK_s\\
&\leq&\int_r^t1_{\{Y_s\geq\bar{L}_s\}}dK_s.
\end{eqnarray*}
\end{proof}

Using Proposition \ref{pr3.1}, we obtain the following comparison results.
\begin{proposition}\label{pr3.2} Let $\bar{g}$ satisfy (\hyperref[2A1]{2A1}). Let $\bar{\xi}\in L_D({\mathcal{F}}_T)$ and $\bar{L}\in{\mathcal{C}}_D$ such that $\bar{\xi}\geq \xi$, $\bar{\xi}\geq\bar{L}_T$, and for all $t\in[0,T]$, $\bar{L}_t\geq L_t.$ Let the BSDE$(g,\xi,L)$ admit a solution $(Y,Z,K)$, and let the BSDE$(g,\xi)$ admit a solution $(y,z)$. Then the following hold:

(i) If the RBSDE$(\bar{g},\bar{\xi},\bar{L}\vee Y)$ admits a solution $(\bar{Y},\bar{Z},\bar{K})$ such that $\bar{g}(t,\bar{Y}_t,\bar{Z}_t)\geq g(t,\bar{Y}_t,\bar{Z}_t)$, $dt\times dP\textrm{-}a.e.,$ then $(\bar{Y},\bar{Z},\bar{K})$ is a solution to the RBSDE$(\bar{g},\bar{\xi},\bar{L})$ such that for all $t\in[0,T]$, $\bar{Y}_t\geq Y_t$;

(ii) If the RBSDE$(\bar{g},\bar{\xi},\bar{L}\vee y)$ admits a solution $(\bar{y},\bar{z},\bar{k})$ such that $\bar{g}(t,\bar{y}_t,\bar{z}_t)\geq g(t,\bar{y}_t,\bar{z}_t)$, $dt\times dP\textrm{-}a.e.,$ then $(\bar{y},\bar{z},\bar{k})$ is a solution to the RBSDE$(\bar{g},\bar{\xi},\bar{L})$ such that for all $t\in[0,T]$, $\bar{y}_t\geq y_t$.
\end{proposition}
\begin{proof} Proof of (i): Since $\bar{Y}_t\geq\bar{L}_t\vee Y_t\geq Y_t$, by Proposition \ref{pr3.1}, we have for any $0\leq r\leq t\leq T,$
\begin{equation*}
\int_r^t1_{\{Y_s\geq \bar{L}_s\}}d\bar{K}_s=\int_r^t1_{\{Y_s\geq \bar{L}_s\vee Y_s\}}d\bar{K}_s\leq\int_r^t1_{\{Y_s\geq \bar{L}_s\vee Y_s\}}d{K}_s=\int_r^t1_{\{Y_s\geq \bar{L}_s\}}dK_s.
\end{equation*}
This together with the assumption that $\bar{L}_t\geq L_t$ and the fact that $\int_0^T1_{\{{Y}_t>L_t\}}dK_t=0$, implies
\begin{eqnarray*}
\int_0^T(\bar{Y}_t-\bar{L}_t)d\bar{K}_t&=&\int_0^T(\bar{Y}_t-(\bar{L}_t\vee{Y}_t))d\bar{K}_t+\int_0^T((\bar{L}_t\vee{Y}_t)-\bar{L}_t)d\bar{K}_t\\
&=&\int_0^T1_{\{{Y}_t\geq \bar{L}_t\}}((\bar{L}_t\vee{Y}_t)-\bar{L}_t)d\bar{K}_t\\
&\leq&\int_0^T1_{\{{Y}_t\geq \bar{L}_t\}}((\bar{L}_t\vee{Y}_t)-\bar{L}_t)dK_t\\
&=&\int_0^T1_{\{{Y}_t>\bar{L}_t\}}((\bar{L}_t\vee{Y}_t)-\bar{L}_t)dK_t\\
&\leq&\int_0^T1_{\{{Y}_t>L_t\}}((\bar{L}_t\vee{Y}_t)-\bar{L}_t)dK_t\\
&=&0.
\end{eqnarray*}
This implies that $(\bar{Y},\bar{Z},\bar{K})$ is a solution to the RBSDE$(\bar{g},\bar{\xi},\bar{L})$. We obtain (i).

Proof of (ii): Since the RBSDE$(g,\xi,y\wedge\bar{L})$ admits a solution $(y,z,0)$, by (i), we obtain (ii).
\end{proof}

To conveniently treat various spaces used to characterize the uniqueness (maximality or minimality) of the solution, we introduce the space ${\mathcal{Y}}_D$, which is a subset of ${\mathcal{C}}_D$ such that:
for any $x^1, x^2\in{\mathcal{Y}}_D$ satisfying $x_t^1\leq x_t^2$ for any $t\in[0,T]$, the set $\{y_t\in{\mathcal{C}}_D:\forall t\in[0,T], x_t^1\leq y_t\leq x_t^2\}$ is contained in ${\mathcal{Y}}_D$. Clearly, ${\mathcal{C}}_D, {\mathcal{S}}_D^p$ and ${\mathcal{S}}_D^\infty$ are some examples of ${\mathcal{Y}}_D$.

Using Proposition \ref{pr3.2}, and Lemmas \ref{l2.3} and \ref{l2.4}, we obtain Proposition \ref{pr3.3}, where the existence of maximal solutions of RBSDEs (resp. BSDEs) plays a crucial role in the study of the uniqueness of solutions in Sections \ref{sec4} and \ref{sec5}.

\begin{proposition}\label{pr3.3} Let $\bar{g}$ satisfy (\hyperref[2A1]{2A1}). Let $\bar{\xi}\in L_D({\mathcal{F}}_T)$ and $\bar{L}, L\in{\mathcal{Y}}_D$ such that $\bar{\xi}\geq \xi$, $\bar{\xi}\geq\bar{L}_T$, and for all $t\in[0,T]$, $\bar{L}_t\geq L_t$. Let the RBSDE$(\bar{g},\bar{\xi},\bar{L})$ admit a maximal solution $(\bar{Y},\bar{Z},\bar{K})$ such that $\bar{Y}\in{\mathcal{Y}}_D$. Assume that

1) $dt\times dP\textrm{-}a.e.$, $\forall(y,z)\in D\times\mathbf{R}$, $\bar{g}(t,\bar{L}_t\vee y,z)\geq g(t,\bar{L}_t\vee y,z)$;

2) $\forall l\in{\mathcal{Y}}_D$ with $l_T\leq\bar{\xi}$, the RBSDE$(\bar{g},\bar{\xi},\bar{L}\vee l)$ admits a solution $(\bar{y}^l,\bar{z}^l,\bar{k}^l)$ such that $\bar{y}^l\in{\mathcal{Y}}_D$.\\
Then the following hold:

(i) The RBSDE$(g,\xi,L)$ admits a maximal solution $(Y,Z,K)$ such that $Y\in{\mathcal{Y}}_D$. Moreover, for all $t\in[0,T]$, $Y_t\leq\bar{Y}_t$;

(ii) If the BSDE$(g,\xi)$ admits a solution $(\tilde{y},\tilde{z})$ such that $\tilde{y}\in{\mathcal{Y}}_D$, then it admits a maximal solution $(\hat{Y},\hat{Z})$ such that $\hat{Y}\in{\mathcal{Y}}_D$. Moreover, for all $t\in[0,T]$, $\hat{Y}_t\leq\bar{Y}_t$.

\begin{proof}
Proof of (i): It is clear that the RBSDE$(g,\xi,L)$ is dominated by $\bar{Y}.$ Then, by Lemma \ref{l2.3}, the RBSDE$(g,\xi,L)$ admits a maximal solution $(Y,Z,K)$ such that $Y_t\leq\bar{Y}_t$ and $Y\in{\mathcal{Y}}_D.$

Since for any solution $(\zeta,\hat{z},\hat{k})$ to the RBSDE$(g,\xi,L)$ such that $\zeta\in{\mathcal{Y}}_D$, the RBSDE$(\bar{g},\bar{\xi},\bar{L}\vee {\zeta})$ admits a solution $(\bar{y}^{\zeta},\bar{z}^{\zeta},\bar{k}^{\zeta})$ such that $\bar{y}^{\zeta}\in{\mathcal{Y}}_D$ and $\bar{g}(t,\bar{y}^{\zeta}_t,\bar{z}^{\zeta}_t)\geq g(t,\bar{y}^{\zeta}_t,\bar{z}^{\zeta}_t)$, $dt\times dP\textrm{-}a.e.$, it follows from Proposition \ref{pr3.2}(i) that $(\bar{y}^\zeta,\bar{z}^\zeta,\bar{k}^{\zeta})$ is a solution to the RBSDE$(\bar{g},\bar{\xi},\bar{L})$ such that $\bar{y}^{\zeta}_t\geq{\zeta}_t$. Hence, we have ${\zeta}_t\leq\bar{Y}_t.$ This implies ${\zeta}_t\leq Y_t.$ We obtain (i).

Proof of (ii): Since for any solution $({\zeta},\hat{z})$ to the BSDE$(g,\xi)$ such that ${\zeta}\in{\mathcal{Y}}_D$, the RBSDE$(\bar{g},\bar{\xi},\bar{L}\vee {\zeta})$ admits a solution $(\bar{y}^{\zeta},\bar{z}^{\zeta},\bar{k}^{\zeta})$ such that $\bar{y}^{\zeta}\in{\mathcal{Y}}_D$ and $\bar{g}(t,\bar{y}^{\zeta}_t,\bar{z}^{\zeta}_t)\geq g(t,\bar{y}^{\zeta}_t,\bar{z}^{\zeta}_t)$, $dt\times dP\textrm{-}a.e.$, it follows from Proposition \ref{pr3.2}(ii) that $(\bar{y}^{\zeta},\bar{z}^{\zeta},\bar{k}^{\zeta})$ is a solution to the RBSDE$(\bar{g},\bar{\xi},\bar{L})$ such that $\bar{y}^{\zeta}_t\geq{\zeta}_t$, and hence, ${\zeta}_t\leq\bar{Y}_t.$

Moreover, since the BSDE$(g,\xi)$ is dominated by ${\zeta}$ and $\bar{Y},$ from Lemma \ref{l2.4}, it follows that the BSDE$(g,\xi)$ admits a maximal solution $(\hat{Y},\hat{Z})$ such that $\hat{Y}_t\leq \bar{Y}_t$ and $\hat{Y}\in{\mathcal{Y}}_D.$ This implies $\zeta_t\leq \hat{Y}_t$. We obtain (ii).
\end{proof}
\end{proposition}

Using Lemmas \ref{l2.3} and \ref{l2.4}, we also obtain the following comparison results.

\begin{proposition}\label{pr3.4} Let $Y^1$ and $Y^2$ be the two semimartingales defined in (\ref{2.4}) and (\ref{2.5}), respectively. Let $Y^1, Y^2, L\in{\mathcal{Y}}_D$ and $Y^2_T\leq\xi\leq Y^1_T$. Let the RBSDE$(g,\xi,L)$ admit a minimal solution $(Y,Z,K)$ such that $Y\in{\mathcal{Y}}_D$, and let the BSDE$(g,\xi)$ admit a maximal solution $(y,z)$ such that $y\in{\mathcal{Y}}_D$. Then the following hold:

(i) If $Y^1_t\geq L_t$ and $h_1(t)\geq g(t,Y_t^1,Z_t^1),\ dt\times dP\textrm{-}a.e$., then for all $t\in[0,T]$, $Y_t^1\geq{Y}_t$;

(ii) If $h_2(t)\leq g(t,Y_t^2,Z_t^2)$, $dt\times dP\textrm{-}a.e.$, and if the RBSDE$(g,\xi,Y^2)$ admits a solution $(\tilde{y},\tilde{z},\tilde{k})$ such that $\tilde{y}\in{\mathcal{Y}}_D$, then for all $t\in[0,T]$, $y_t\geq Y_t^2$.
\end{proposition}
\begin{proof} Proof of (i): Since $Y^1_T\geq\xi$, and $dt\times dP\textrm{-}a.e.$, $Y^1_t\geq L_t$ and $h_1(t)\geq g(t,Y_t^1,Z_t^1)$, we get that the RBSDE$(g,\xi,L)$ is dominated by $Y^1$. Then by Lemma \ref{l2.3}, the RBSDE$(g,\xi,L)$ admits a solution $(\hat{Y},\hat{Z},\hat{K})$ such that $L_t\leq \hat{Y}_t\leq Y_t^1$ and $\hat{Y}\in{\mathcal{Y}}_D.$ This implies (i).

Proof of (ii): Since $Y^2_t\leq\tilde{y}_t$ and $h_2(t)\leq g(t,Y_t^2,Z_t^2),\ dt\times dP\textrm{-}a.e.$, we get that the BSDE$(g,\xi)$ is dominated by $Y^2$ and $\tilde{y}$. Then by Lemma \ref{l2.4}, the BSDE$(g,\xi)$ admits a solution $(\hat{y},\hat{z})$ such that $Y^2_t\leq \hat{y}_t\leq \tilde{y}_t$ and $\hat{y}\in{\mathcal{Y}}_D.$ This implies (ii).
\end{proof}

\begin{remark}\label{r3.5} Propositions \ref{pr3.3} and \ref{pr3.4} mainly rely on the existence, maximality and minimality of solutions. These roughly imply that the comparison theorem for a quadratic RBSDE (resp. BSDE) holds whenever such an RBSDE (resp. BSDE) admits a unique solution in a certain space. Using a different method, this phenomenon was previously observed in \cite[Theorem 2.7]{Zheng1} for BSDEs with generators that are Lipschitz in $y$, in a general setting.
\end{remark}
\section{The existence and uniqueness of bounded solutions}\label{sec4}
\subsection{Case of a one-sided superlinear growth in $y$}
Let ${\mathcal{L}}$ be the set consisting of all the continuous function $l(x):\mathbf{R}\mapsto(0,\infty)$ such that
\begin{equation*}
\int_{-\infty}^0\frac{1}{l(x)}dx=\int_0^\infty\frac{1}{l(x)}dx=\infty.
\end{equation*}
It was introduced by Lepeltier and San Martin \cite[Lemma 1]{LS}. We consider the following assumptions for $g$:
\begin{itemize}
  \item \textbf{(4A1)}\label{4A1} There exist three functions $u(t)\in L^1_+[0,T]$, $l(x)\in{\mathcal{L}}$, $f(x)\in C_+(D),$ and two constants $c\geq b>0$, such that $dt\times dP\textrm{-}a.e.,$ for each $(y,z)\in D\times{\mathbf{R}}^{\mathit{d}},$
\begin{itemize}
  \item\textbf{(i)}\ \ \ \ \ \ \ \ \ \ \ \ \ \ \ \ \ \ \ \ \ \ \ \ \ \ \ \ \ \ \ \ \ \ $1_{\{y\geq c\}}{g}(t,y,0)\leq u(t)l(y);$
  \item\textbf{(ii)} \ \ \ \ \ \ \ \ \ \ \ \ \ \ \ \ \ \ \ \ \ \ \ \ \ \ \ \ \ \ $1_{\{y\leq -c\}}{g}(t,y,0)\geq -u(t)l(y);$
  \item\textbf{(iii)}\ \ \ \ \ \ \ \ \ \ \ \ \ \ \ \ \ \ \ \ \ \ \ \ \ \ \ \ \ $1_{\{y\leq b\}}{g}(t,y,0)\geq -u(t)|y\ln(|y|)|;$
  \item\textbf{(iv)}\ \ \ \ \ \ \ \ \ \ \ \ \ \ \ \ \ \ \ \ \ \ \ \ \ \ $1_{\{y\geq c\}}{g}(t,y,z)\leq u(t)l(y)+f(y)|z|^2.$
\end{itemize}
  \item \textbf{(4A2)}\label{4A2} $g$ satisfies the $\theta$-domination condition for ${\mathcal{S}}_D^\infty$ and $H^2_d$ with $G$, and there exist $\tilde{u}(t)\in L^1_+[0,T]$, $\tilde{l}(x)\in{\mathcal{L}}$, $\tilde{f}(x)\in C_+(D),$ $\tilde{h}(x,y)\in C_+(D\times D),$ and a constant $\tilde{c}>0$, such that $dt\times dP\textrm{-}a.e.,$ for any $y^1,y^2\in{\mathcal{S}}_D^\infty$ and $(y,z)\in D\times{\mathbf{R}}^{\mathit{d}},$
      \begin{equation*}
      1_{\{y\geq \tilde{c}\}}G(t,y,z,y^1_t,y^2_t)\leq \tilde{h}(y^1_t,y^2_t)\tilde{u}(t)\tilde{l}(y)+\tilde{f}(y)|z|^2.
      \end{equation*}
  \item \textbf{(4A2')}\label{4A2'} $dt\times dP\textrm{-}a.e.,$ $g(t,\cdot,\cdot)$ is convex on $D\times{\mathbf{R}}^{\mathit{d}}$.
  \item \textbf{(4A2'')}\label{4A2''} $dt\times dP\textrm{-}a.e.,$ for each $y\in D$, $g(t,y,\cdot)$ is convex in $z$, and there exists a function $h(x,y)\in C_+(D\times D)$, such that for each $y_1, y_2\in D$ and $z\in\mathbf{R}^d$,
      $$|{g}(t,y_1,z)-{g}(t,y_2,z)|\leq h(y_1,y_2)|y_1-y_2|.$$
\end{itemize}

Intuitively, (\hyperref[4A2]{4A2}) is not very easy to be checked. We show some typical cases of (\hyperref[4A2]{4A2}), which are inspired by \cite{FH}.
\begin{remark}\label{r4.1}
For two processes $y^1,y^2\in{\mathcal{C}}_D$, two processes $z^1,z^2\in H^2_d$ and a constant $\theta\in(0,1)$, we set $\Delta_\theta y_t:=\frac{y_t^1-\theta y_t^2}{1-\theta}$ and $\Delta_\theta z_t:=\frac{z_t^1-\theta z_t^2}{1-\theta}$.
\begin{itemize}
\item (i) (\hyperref[4A2']{4A2'}) combined with (\hyperref[4A1]{4A1}-(iv)) implies (\hyperref[4A2]{4A2}). In fact, when (\hyperref[4A2']{4A2'}) holds and the range of $\Delta_\theta y_t$ is included in $D$, we have
\begin{align*}
g(t,y_t^1,z_t^1)-\theta g(t,y_t^2,z_t^2)&=g(t,\theta y_t^2+(1-\theta)\Delta_\theta y_t,\theta z_t^2+(1-\theta)\Delta_\theta z_t)-\theta g(t,y_t^2,z_t^2)\\
&\leq(1-\theta)g(t,\Delta_\theta y_t,\Delta_\theta z_t).
\end{align*}
We choose $G(t,y,z,y_t^1,y_t^2)=g(t,y,z).$ Thus, if (\hyperref[4A1]{4A1}-(iv)) also holds, then (\hyperref[4A2]{4A2}) holds.

\item (ii) (\hyperref[4A2'']{4A2''}) combined with (\hyperref[4A1]{4A1}-(iv)) implies (\hyperref[4A2]{4A2}). In fact, when (\hyperref[4A2'']{4A2''}) holds, we get that for each $\beta>0$,
\begin{align*}
  &g(t,y_t^1,z_t^1)-\theta g(t,y_t^2,z_t^2)\\ &\leq g(t,y_t^1,z_t^1)-g(t,y_t^2,z_t^1)+g(t,y_t^2,z_t^1)-\theta g(t,y_t^2,z_t^2)\\
  &\leq h(y_t^1,y_t^2)|y^1_t-y^2_t|+g(t,y_t^2,\theta z_t^2+(1-\theta)\Delta_\theta z_t)-\theta g(t,y_t^2,z_t^2)\\
  &\leq h(y_t^1,y_t^2)|y^1_t-\theta y^2_t|+(1-\theta)(h(y_t^1,y_t^2)|y_t^2|+h(y_t^2,\beta)(|y_t^2|+\beta)+g(t,\beta,\Delta_\theta z_t))\\
  &\leq (1-\theta)(h(y_t^1,y_t^2)|y_t^2|+h(y_t^2,\beta)(|y_t^2|+\beta)+h(y_t^1,y_t^2)|\Delta_\theta y_t|+g(t,\beta,\Delta_\theta z_t)),\tag{4.1}\label{4.1}
\end{align*}
We choose
\begin{equation*}
G(t,y,z,y_t^1,y_t^2)=h(y_t^1,y_t^2)|y_t^2|+h(y_t^2,\beta)(|y_t^2|+\beta)+h(y_t^1,y_t^2)|y|+g(t,\beta,z).
\end{equation*}
If (\hyperref[4A1]{4A1}-(iv)) also holds and we choose $\beta\geq c$, then we have
\begin{equation*}
G(t,y,z,y^1_t,y^2_t)\leq \bar{h}(y_t^1,y^2_t)\bar{u}(t)(|y|+1)+f(\beta)|z|^2,
\end{equation*}
where $\bar{u}(t)\in L^1_+[0,T]$ is dependent only on $u(t)$, and $\bar{h}\in C_+(D\times D)$ is dependent only on $h$, $l$ and $\beta$. Hence, (\hyperref[4A2]{4A2}) holds. Moreover, if the range of $\Delta_\theta y_t$ is included in $D$, then replacing $\beta$ in (\ref{4.1}) by $\Delta_\theta y_t$, we also have
\begin{align*}
 &g(t,y_t^1,z_t^1)-\theta g(t,y_t^2,z_t^2)\\
  &\leq (1-\theta)(h(y_t^1,y_t^2)|y_t^2|+h(y_t^2,\Delta_\theta y_t)(|y_t^2|+|\Delta_\theta y_t|)+h(y_t^1,y_t^2)|\Delta_\theta y_t|+g(t,\Delta_\theta y_t,\Delta_\theta z_t)),
\end{align*}
which implies that if we assume that the function $h(x,y)$ is a nonnegative process $h\in\mathcal{C}$ independent of $(x,y)$ (i.e., (\hyperref[5A2']{5A2'}) in Section \ref{sec5} holds), then (\hyperref[5A2]{5A2}) in Section \ref{sec5} holds, since we can choose
\begin{equation*}
G(t,y,z,y^1_t,y^2_t)=2h_t|y_t^2|+2h_t|y|+g(t,y,z).\tag{4.2}\label{4.2}
\end{equation*}

\item (iii) Let $g_1(z):{\bf{R}}^d\rightarrow{\bf{R}}$ be a bounded Lipschitz function with a bounded support. From the proof of \cite[Proposition 3.5(iii)]{FH}, we get that there exists a constant $M>0$ such that $g_1(z_t^1)-\theta g_1(z_t^2)\leq(1-\theta)M(1+|\Delta_\theta z_t|).$ Hence, by (ii), if $g$ satisfies (\hyperref[4A2'']{4A2''}) and (\hyperref[4A1]{4A1}-(iv)), then $g(t,y,z)+g_1(z)$ satisfies (\hyperref[4A2]{4A2}). Clearly, $g+g_1$ may be not convex.
\end{itemize}
\end{remark}

We have the following well-posedness results for bounded solutions of RBSDEs.

\begin{theorem} \label{th4.2} Let (\hyperref[4A1]{4A1}-(i)) hold, $\xi\in L_D^\infty({\mathcal{F}}_T)$, and $L\in{\mathcal{S}}_D^\infty$.
Then the RBSDE$(g,\xi,L)$ admits a minimal solution $(\underline{Y},\underline{Z},\underline{K})$ such that $\underline{Y}\in{\mathcal{C}}_D$. Moreover, the following hold:

(i) If (\hyperref[4A1]{4A1}-(iv)) further holds, then the RBSDE$(g,\xi,L)$ admits a maximal solution $(\overline{Y},\overline{Z},\overline{K})$ such that $\overline{Y}\in{\mathcal{S}}_D^\infty$;

(ii) If (\hyperref[4A2]{4A2}) further holds, then the RBSDE$(g,\xi,L)$ admits a unique solution $(Y,Z,K)$ such that $Y\in{\mathcal{S}}_D^\infty$.
\end{theorem}
\begin{proof} We divide this proof into three steps.

\textbf{Step 1.} We assume that the ranges of $\xi$ and $L$ are both included in $[b_1,c_1]\subset D$ and that $c_2\geq c\vee c_1.$ By Fan \cite[Lemma 3.1]{F}, the backward ODE:
\begin{equation*}
\varphi(t)=c_2+\int_t^Tu(s)l(\varphi(s))ds, \ \ t\in[0,T]
\end{equation*}
admits a unique solution $\varphi(t)$ such that for each $t\in[0,T]$, $c_2\leq\varphi(t)\leq\varphi(0).$ By (\hyperref[4A1]{4A1}-(i)), we have
\begin{equation*}
g(t,\varphi(t),0)\leq u(t)l(\varphi(t)),\ \ \ dt\times dP\textrm{-}a.e.,
\end{equation*}
which together with the facts that $\xi\leq c_2$ and $L_t\leq\varphi(t)$, implies that the RBSDE$(g,\xi,L)$ is dominated by $\varphi(t).$ It follows from Lemma \ref{l2.3} that the RBSDE$(g,\xi,L)$ admits a minimal solution $(\underline{Y},\underline{Z},\underline{K})$ such that $\underline{Y}\in{\mathcal{C}}_D$. Moreover, for each $t\in[0,T],$ $b_1\leq \underline{Y}_t\leq \varphi(0).$

\textbf{Step 2.} Proof of (i): By \cite[Theorem 3.1]{F}, the BSDE$(u(t)l(y)+f(y)|z|^2,c_2)$ admits a maximal solution $(\varphi(t),0)\in{\mathcal{S}}_D^\infty\times{\mathcal{H}}^{BMO}_d$. Lemma \ref{lmB.1} in Appendix B further implies that $(\varphi(t),0)$ is a maximal solution to the BSDE$(u(t)l(y)+f(y)|z|^2,c_2)$ such that $\varphi(t)\in{\mathcal{S}}_D^\infty$. Since $\varphi(t)\geq c_2\geq L_t\vee c$, we get that $(\varphi(t),0,0)$ is a solution to the RBSDE$(u(t)l(y)+f(y)|z|^2,c_2,L\vee c)$ such that $\varphi(t)\in{\mathcal{S}}_D^\infty.$

Let $(\hat{y},\hat{z},\hat{k})$ be another solution to the RBSDE$(u(t)l(y)+f(y)|z|^2,c_2,L\vee c)$ such that $\hat{y}\in{\mathcal{S}}_D^\infty.$ Set
\begin{equation*}
M:=\varphi(0)\vee\sup_{t\in[0,T]}\|\hat{y}_t\|_\infty.
\end{equation*}
By \cite[Theorem 3.1]{F} again, the BSDE$(u(t)l(y)+f(y)|z|^2,M)$ admits a maximal solution $(\varphi^M(t),0)\in{\mathcal{S}}_D^\infty\times{\mathcal{H}}^{BMO}_d$, and moreover $\varphi^M(t)\geq \varphi(t)$ and $\varphi^M(t)\geq \hat{y}_t$. It follows that the RBSDE$(u(t)l(y)+f(y)|z|^2,c_2,L\vee c)$ is dominated by $\varphi^M(t)$. Lemma \ref{l2.3} then implies that the RBSDE$(u(t)l(y)+f(y)|z|^2,c_2,L\vee c)$ admits a maximal solution $(\bar{y},\bar{z},\bar{k})$ such that $\bar{y}\leq\varphi^M(t)$, which gives that $\bar{y}_t\geq \varphi(t)$ and $\bar{y}_t\geq \hat{y}_t$. By Proposition \ref{pr3.1}, we obtain that $\bar{k}_t\leq0$, and thus $(\bar{y},\bar{z})$ is a solution to the BSDE$(u(t)l(y)+f(y)|z|^2,c_2)$ such that $\bar{y}\in{\mathcal{S}}_D^\infty.$ This implies $\varphi(t)\geq \bar{y}_t\geq \hat{y}_t$. Hence, $(\varphi(t),0,0)$ is a maximal solution to the RBSDE$(u(t)l(y)+f(y)|z|^2,c_2,L\vee c)$ such that $\varphi(t)\in{\mathcal{S}}_D^\infty.$

Moreover, for each $l\in{\mathcal{S}}^\infty_D$ such that $l_T\leq c_2$, the consequence of Step 1 implies that the RBSDE$(u(t)l(y)+f(y)|z|^2,c_2,L\vee c\vee l)$ admits a solution $(y^l,z^l,k^l)$ such that $y^l\in{\mathcal{S}}^\infty_D$. Then, by (\hyperref[4A1]{4A1}-(iv)) and Proposition \ref{pr3.3}(i) with ${\mathcal{Y}}_D={\mathcal{S}}^\infty_D$, we obtain (i).

\textbf{Step 3.}  Proof of (ii): By Step 1, the RBSDE$(g,\xi,L)$ admits a minimal solution $(\underline{Y},\underline{Z},\underline{K})$ such that ${\underline{Y}\in\mathcal{S}}_D^\infty$. Let $(\hat{Y},\hat{Z},\hat{K})$ be another solution to the RBSDE$(g,\xi,L)$ such that ${\hat{Y}\in\mathcal{S}}_D^\infty$. Note that this proof requires the fact that $\hat{Y}_t\geq \underline{Y}_t$. Since $\hat{Y}_t\geq \underline{Y}_t$, for each $\theta\in(0,1)$, we have
\begin{equation}\label{4.3}
\frac{\hat{Y}-\theta \underline{Y}}{1-\theta}=\hat{Y}+\frac{\theta}{1-\theta}(\hat{Y}-\underline{Y})\in{\mathcal{C}}_D.\tag{4.3}
\end{equation}
For $\theta\in(0,1)$, we have
\begin{eqnarray*}
\frac{\hat{Y}_t-\theta \underline{Y}_t}{1-\theta}&=&\xi+\int_t^T\frac{\hat{g}_s}{1-\theta}ds+\frac{(\hat{K}_T-\theta \underline{K}_T)-(\hat{K}_t-\theta \underline{K}_t)}{1-\theta}-\int_t^T\frac{\hat{Z}_s-\theta \underline{Z}_s}{1-\theta}\cdot dB_s
\end{eqnarray*}
where
\begin{equation*}
\frac{\hat{g}_s}{1-\theta}:=\frac{1}{{1-\theta}}(g(s,\hat{Y}_s,\hat{Z}_s)-\theta g(s,\underline{Y}_s,\underline{Z}_s)).
\end{equation*}
Since $\hat{Y}$ and $\underline{Y}$ both belong to ${\mathcal{S}}_D^\infty$, by (\ref{4.3}) and (\hyperref[4A2]{4A2}), we have
\begin{equation*}\label{4.4}
 \forall\theta\in(0,1), \ \  \frac{\hat{g}_t}{1-\theta}\leq G\left(t,\frac{\hat{Y}_t-\theta \underline{Y}_t}{1-\theta},\frac{\hat{Z}_t-\theta \underline{Z}_t}{1-\theta},\hat{Y}_t,\underline{Y}_t\right),\tag{4.4}
\end{equation*}
and $dt\times dP\textrm{-}a.e.,$ for each $(y,z)\in D\times{\mathbf{R}}^{\mathit{d}},$
\begin{equation*}
1_{\{y\geq \tilde{c}\}}G(t,y,z,\hat{Y}_t,\underline{Y}_t)\leq \tilde{h}(\hat{Y}_t,\underline{Y}_t)\tilde{u}(t)\tilde{l}(y)+\tilde{f}(y)|z|^2.
\end{equation*}
By (\ref{4.3}), we have for each $\theta\in(0,1)$, $\frac{\hat{Y}-\theta \underline{Y}}{1-\theta}\in{\mathcal{S}}_D^\infty$. This together with the conclusion of Step 1 implies that for each $\theta\in(0,1)$, the RBSDE$(G(t,y,z,\hat{Y}_t,\underline{Y}_t),\xi,\frac{\hat{Y}-\theta \underline{Y}}{1-\theta})$ admits a minimal solution $(y^\theta,z^\theta,k^\theta)$ such that $y^\theta\in{\mathcal{S}}_D^\infty.$

In the following, we will show that for each $\theta\in(0,1)$, $(y^\theta,z^\theta,k^\theta)$ is actually a solution to the RBSDE$(G(t,y,z,\hat{Y}_t,\underline{Y}_t),\xi,\hat{Y})$. This is crucial to this proof. Since for each $\theta\in(0,1)$ and $s\in[0,T]$, $(y^\theta_s-\frac{\hat{Y}_s-\theta \underline{Y}_s}{1-\theta})^+=(y^\theta_s-\frac{\hat{Y}_s-\theta \underline{Y}_s}{1-\theta}),$ by applying Tanaka's formula to $(y^\theta_s-\frac{\hat{Y}_s-\theta \underline{Y}_s}{1-\theta})^+-(y^\theta_s-\frac{\hat{Y}_s-\theta \underline{Y}_s}{1-\theta}),$ we deduce that for any $0\leq r\leq t\leq T,$
\begin{align}
0&=\left(y^\theta_t-\frac{\hat{Y}_t-\theta \underline{Y}_t}{1-\theta}\right)^+-\left(y^\theta_t-\frac{\hat{Y}_t-\theta \underline{Y}_t}{1-\theta}\right)\notag\\
&\ \ \ \ \ \ \ \ \ \ -\left(\left(y^\theta_r-\frac{\hat{Y}_r-\theta \underline{Y}_r}{1-\theta}\right)^+-\left(y^\theta_r-\frac{\hat{Y}_r-\theta \underline{Y}_r}{1-\theta}\right)\right)\notag\\
&=\int_r^t1_{\{y^\theta_s=\frac{\hat{Y}_s-\theta \underline{Y}_s}{1-\theta}\}}\left(G(t,y^\theta_s,z^\theta_s,\hat{Y}_s,\underline{Y}_s )-\frac{\hat{g}_s}{1-\theta}\right)ds+\int_r^t1_{\{y^\theta_s=\frac{\hat{Y}_s-\theta \underline{Y}_s}{1-\theta}\}}dk^\theta_s\tag{4.5}\label{4.5}\\
&\ \ \ \ \ -\int_r^t1_{\{y^\theta_s=\frac{\hat{Y}_s-\theta \underline{Y}_s}{1-\theta}\}}d\left(\frac{\hat{K}_s-\theta \underline{K}_s}{1-\theta}\right)-\int_r^t1_{\{y^\theta_s=\frac{\hat{Y}_s-\theta \underline{Y}_s}{1-\theta}\}}\left(z^\theta_s-\frac{\hat{Z}_s-\theta \underline{Z}_s}{1-\theta}\right)\cdot dB_s\notag\\
&\ \ \ \ \ +\frac{1}{2}\ell_t^0\left(y^\theta_s-\frac{\hat{Y}_s-\theta \underline{Y}_s}{1-\theta}\right)-\frac{1}{2}\ell_r^0\left(y^\theta_s-\frac{\hat{Y}_s-\theta \underline{Y}_s}{1-\theta}\right),\notag
\end{align}
where $\ell_t^0\left(y^\theta_s-\frac{\hat{Y}_s-\theta \underline{Y}_s}{1-\theta}\right)$ is the local time of $y^\theta_s-\frac{\hat{Y}_s-\theta \underline{Y}_s}{1-\theta}$ at time $t$ and level $0$. This implies that
\begin{equation*}
\int_r^t1_{\{y^\theta_s=\frac{\hat{Y}_s-\theta \underline{Y}_s}{1-\theta}\}}\left(z^\theta_s-\frac{\hat{Z}_s-\theta \underline{Z}_s}{1-\theta}\right)\cdot dB_s=0,
\end{equation*}
which further gives $dt\times dP\textrm{-}a.e.,$
\begin{equation*}
1_{\{y^\theta_t=\frac{\hat{Y}_t-\theta \underline{Y}_t}{1-\theta}\}}G(t,y^\theta_t,z^\theta_t,\hat{Y}_t,\underline{Y}_t) =1_{\{y^\theta_t=\frac{\hat{Y}_t-\theta \underline{Y}_t}{1-\theta}\}}G\left(t,\frac{\hat{Y}_t-\theta \underline{Y}_t}{1-\theta},\frac{\hat{Z}_t-\theta \underline{Z}_t}{1-\theta},\hat{Y}_t,\underline{Y}_t\right).
\end{equation*}
This together with (\ref{4.4}) and (\ref{4.5}), implies
\begin{equation*}
\int_r^t1_{\{y^\theta_s=\frac{\hat{Y}_s-\theta \underline{Y}_s}{1-\theta}\}}dk^\theta_s\leq\int_r^t1_{\{y^\theta_s=\frac{\hat{Y}_s-\theta \underline{Y}_s}{1-\theta}\}}d\left(\frac{\hat{K}_s-\theta \underline{K}_s}{1-\theta}\right).
\end{equation*}
Since $\int_0^T1_{\{y^\theta_s>\frac{\hat{Y}_s-\theta \underline{Y}_s}{1-\theta}\}}dk^\theta_s=0$, the inequality above implies
\begin{equation*}
k^\theta_t-k^\theta_r=\int_r^t1_{\{y^\theta_s=\frac{\hat{Y}_s-\theta \underline{Y}_s}{1-\theta}\}}dk^\theta_s+\int_r^t1_{\{y^\theta_s>\frac{\hat{Y}_s-\theta \underline{Y}_s}{1-\theta}\}}dk^\theta_s\leq\int_r^t1_{\{y^\theta_s=\frac{\hat{Y}_s-\theta \underline{Y}_s}{1-\theta}\}}d\left(\frac{\hat{K}_s-\theta \underline{K}_s}{1-\theta}\right).
\end{equation*}
From this and the fact that $\int_r^t1_{\{\hat{Y}_s>L_s\}}d\hat{K}_s=0$, it follows that
\begin{align}
\int_r^t1_{\{\hat{Y}_s>\underline{Y}_s\}}dk^\theta_s&\leq\int_r^t1_{\{\hat{Y}_s>\underline{Y}_s\}}1_{\{y^\theta_s=\frac{\hat{Y}_s-\theta \underline{Y}_s}{1-\theta}\}}d\left(\frac{\hat{K}_s-\theta \underline{K}_s}{1-\theta}\right)\notag\\
&=\frac{1}{1-\theta}\left[\int_r^t1_{\{\hat{Y}_s>\underline{Y}_s\}}1_{\{y^\theta_s=\frac{\hat{Y}_s-\theta \underline{Y}_s}{1-\theta}\}}d\hat{K}_s-\theta\int_r^t1_{\{\hat{Y}_s>\underline{Y}_s\}}1_{\{y^\theta_s=\frac{\hat{Y}_s-\theta \underline{Y}_s}{1-\theta}\}}d\underline{K}_s\right]\notag\\
&\leq\frac{1}{1-\theta}\left[\int_r^t1_{\{\hat{Y}_s>L_s\}}d\hat{K}_s-\theta\int_r^t1_{\{\hat{Y}_s>\underline{Y}_s\}}1_{\{y^\theta_s=\frac{\hat{Y}_s-\theta \underline{Y}_s}{1-\theta}\}}d\underline{K}_s\right]\notag\\
&=-\frac{\theta}{1-\theta}\int_r^t1_{\{\hat{Y}_s>\underline{Y}_s\}}1_{\{y^\theta_s=\frac{\hat{Y}_s-\theta \underline{Y}_s}{1-\theta}\}}d\underline{K}_s\notag\\
&\leq0,\notag
\end{align}
which implies $\int_r^t1_{\{\hat{Y}_s>\underline{Y}_s\}}dk^\theta_s=0.$ Then, by the facts that $\hat{Y}_t\geq \underline{Y}_t$ and $\int_0^T1_{\{y^\theta_s>\frac{\hat{Y}_s-\theta \underline{Y}_s}{1-\theta}\}}dk^\theta_s=0$, we have for any $0\leq r\leq t\leq T,$
\begin{eqnarray*}
\int_r^t1_{\{\hat{Y}_s=\underline{Y}_s\}}dk^\theta_s&=&\int_r^t1_{\{\hat{Y}_s=\underline{Y}_s\}}dk^\theta_s+\int_r^t1_{\{\hat{Y}_s>\underline{Y}_s\}}dk^\theta_s\\
&=&k^\theta_t-k^\theta_r\\
&=&\int_r^t1_{\{y^\theta_s=\frac{\hat{Y}_s-\theta \underline{Y}_s}{1-\theta}\}}dk^\theta_s.
\end{eqnarray*}
It follows that $1_{\{\hat{Y}_t=\underline{Y}_t\}}=1_{\{y^\theta_t=\frac{\hat{Y}_t-\theta \underline{Y}_t}{1-\theta}\}},\ dk^\theta_s\times dP\textrm{-}a.e.,$ and thus
\begin{eqnarray*}
\int_0^T(y^\theta_s-\hat{Y}_s)1_{\{y^\theta_s=\frac{\hat{Y}_s-\theta \underline{Y}_s}{1-\theta}\}}dk^\theta_s
&=&\int_0^T(y^\theta_s-\hat{Y}_s)1_{\{\hat{Y}_s=\underline{Y}_s\}}1_{\{y^\theta_s=\frac{\hat{Y}_s-\theta \underline{Y}_s}{1-\theta}\}}dk^\theta_s\\
&=&\int_0^T(y^\theta_s-\hat{Y}_s)1_{\{y^\theta_s=\hat{Y}_s\}}1_{\{y^\theta_s=\frac{\hat{Y}_s-\theta \underline{Y}_s}{1-\theta}\}}dk^\theta_s\\
&=&0.
\end{eqnarray*}
This together with the fact that $\int_0^T1_{\{y^\theta_s>\frac{\hat{Y}_s-\theta \underline{Y}_s}{1-\theta}\}}dk^\theta_s=0$ gives
\begin{align*}
\int_0^T(y^\theta_s-\hat{Y}_s)dk^\theta_s=\int_0^T(y^\theta_s-\hat{Y}_s)1_{\{y^\theta_s=\frac{\hat{Y}_s-\theta \underline{Y}_s}{1-\theta}\}}dk^\theta_s+\int_0^T(y^\theta_s-\hat{Y}_s)1_{\{y^\theta_s>\frac{\hat{Y}_s-\theta \underline{Y}_s}{1-\theta}\}}dk^\theta_s=0\tag{4.6}\label{4.6}
\end{align*}
Moreover, since $\hat{Y}_t\geq \underline{Y}_t,$ we have for each $t\in[0,T]$,
\begin{equation*}
y^\theta_t\geq\frac{\hat{Y}_t-\theta \underline{Y}_t}{1-\theta}=\hat{Y}_t+\frac{\theta}{1-\theta}(\hat{Y}_t-\underline{Y}_t)\geq \hat{Y}_t.
\end{equation*}
Then, from this and (\ref{4.6}), we get that for each $\theta\in(0,1)$, $(y^\theta,z^\theta,k^\theta)$ is actually a solution to the RBSDE$(G(t,y,z,\hat{Y}_t,\underline{Y}_t),\xi,\hat{Y})$. By the conclusion of Step 2, the RBSDE$(G(t,y,z,\hat{Y}_t,\underline{Y}_t),\xi,\hat{Y})$ admits a maximal solution $(\tilde{y},\tilde{z},\tilde{k})$ such that $\tilde{y}\in{\mathcal{S}}_D^\infty$. This implies that for each $\theta\in(0,1)$, $(1-\theta)\tilde{y}_t\geq(1-\theta)y^\theta_t\geq \hat{Y}_t-\theta \underline{Y}_t$. When $\theta$ tends to $1$, we have $\underline{Y}_t\geq \hat{Y}_t,$ which implies $\hat{Y}_t=\underline{Y}_t.$ From this, we further derive that $(\hat{Y}_t,\hat{Z}_t,\hat{K}_t)=(\underline{Y}_t,\underline{Z}_t,\underline{K}_t),\ dt\times dP\textrm{-}a.e.$ We obtain (ii).
\end{proof}

We have the following well-posedness results for bounded solutions of BSDEs.
\begin{theorem}\label{th4.3}
Let one of the following two conditions hold:

(i) $D=\mathbf{R}$, $\xi\in L_D^\infty({\mathcal{F}}_T)$, and (\hyperref[4A1]{4A1}-(i)(ii)) hold;

(ii) $D=(0,\infty)$, $\xi\in L_D^\infty({\mathcal{F}}_T)$, and (\hyperref[4A1]{4A1}-(i)(iii)) hold.\\
Then the BSDE$(g,\xi)$ admits at least one solution $(Y,Z)$ such that ${Y\in\mathcal{S}}_D^\infty$. Moreover, if (\hyperref[4A2]{4A2}) also holds, then the BSDE$(g,\xi)$ admits a unique solution $(Y,Z)$ such that ${Y\in\mathcal{S}}_D^\infty$.
\end{theorem}
\begin{proof} We assume that the range of $\xi$ is included in a closed subset $[b_1,c_1]$ of $D$. Let $b_2, b_3, c_2$ be the constants such that $b_2\leq b_1\wedge (-c)$, $0<b_3<b_1\wedge b\wedge1$ and $c_2\geq c_1\vee c$.

\textbf{Case (i):} Let (\hyperref[4A1]{4A1}-(i)(ii)) hold and $D=\mathbf{R}$. By \cite[Lemma 3.1]{F}, the following two backward ODEs:
\begin{equation*}
\phi(t)=b_2-\int_t^Tu(s)l(\phi(s))ds\ \ \ \textrm{and}\ \ \ \varphi(t)=c_2+\int_t^Tu(s)l(\varphi(s))ds, \ \ t\in[0,T]
\end{equation*}
have unique solutions $\phi(t)$ and $\varphi(t)$ respectively, such that for each $t\in[0,T]$,
\begin{equation}\label{4.7}
\phi(0)\leq\phi(t)\leq b_2\leq c_2\leq\varphi(t)\leq\varphi(0).\tag{4.7}
\end{equation}
Then by (\hyperref[4A1]{4A1}-(i)(ii)), we have
\begin{equation*}
g(t,\varphi(t),0)\leq u(t)l(\varphi(t)) \ \ \textrm{and} \ \ g(t,\phi(t),0)\geq -u(t)l(\phi(t)),\ \ \ dt\times dP\textrm{-}a.e.,
\end{equation*}
which together with the facts that $b_2\leq\xi\leq c_2$ and $\phi(t)\leq\varphi(t)$, implies that the BSDE$(g,\xi)$ is dominated by $\phi(t)$ and $\varphi(t)$. It follows from Lemma \ref{l2.4} that the BSDE$(g,\xi)$ admits a solution $(Y_t,Z_t)$ such that $\phi(0)\leq Y_t\leq\varphi(0).$

Let us further assume that (\hyperref[4A2]{4A2}) holds. Let $(\hat{Y},\hat{Z})$ be another solution to the BSDE$(g,\xi)$ such that ${\hat{Y}\in\mathcal{S}}^\infty$. For $\theta\in(0,1)$, we have
\begin{equation}\label{4.8}
\frac{\hat{Y}_t-\theta Y_t}{1-\theta}=\xi+\int_t^T\frac{1}{{1-\theta}}(g(s,\hat{Y}_s,\hat{Z}_s)-\theta g(s,Y_s,Z_s))ds-\int_t^T\frac{\hat{Z}_s-\theta Z_s}{1-\theta}\cdot dB_s.\tag{4.8}
\end{equation}
Since $D=\mathbf{R}$, and $\hat{Y}_t$ and $Y_t$ both belong to ${\mathcal{S}}^\infty$, by (\hyperref[4A2]{4A2}), we have
\begin{equation}\label{4.9}
\frac{1}{{1-\theta}}(g(t,\hat{Y}_t,\hat{Z}_t)-\theta g(t,Y_t,Z_t))\leq G\left(t,\frac{\hat{Y}_t-\theta Y_t}{1-\theta},\frac{\hat{Z}_t-\theta Z_t}{1-\theta},\hat{Y}_t,Y_t\right),\tag{4.9}
\end{equation}
and $dt\times dP\textrm{-}a.e.,$ for each $(y,z)\in D\times{\mathbf{R}}^{\mathit{d}},$
\begin{equation*}
1_{\{y\geq \tilde{c}\}}G(t,y,z,\hat{Y}_t,Y_t)\leq \tilde{h}(\hat{Y}_t,Y_t)\tilde{u}(t)\tilde{l}(y)+\tilde{f}(y)|z|^2.
\end{equation*}
By Theorem \ref{th4.2}, for each $\theta\in(0,1)$, the RBSDE$(G(t,y,z,\hat{Y}_t,Y_t),\xi,\frac{\hat{Y}-\theta Y}{1-\theta})$ admits a minimal solution $(y^\theta,z^\theta,k^\theta)$ such that $y^\theta\in{\mathcal{S}}^\infty$.

In view of $\frac{\hat{Y}_t-\theta Y_t}{1-\theta}\leq y^\theta_t$, (\ref{4.8}) and (\ref{4.9}), we get that the BSDE$(G(t,y,z,\hat{Y}_t,Y_t),\xi)$ is dominated by $\frac{\hat{Y}-\theta Y}{1-\theta}$ and $y^\theta$. Lemma \ref{l2.4} then implies that the BSDE$(G(t,y,z,\hat{Y}_t,Y_t),\xi)$ admits a solution $(\hat{y}^\theta,\hat{z}^\theta)$ such that $\hat{y}^\theta_t\geq\frac{\hat{Y}_t-\theta Y_t}{1-\theta}$ and $\hat{y}^\theta\in{\mathcal{S}}^\infty$. By Theorem \ref{th4.2}(i), the RBSDE$(G(t,y,z,\hat{Y}_t,Y_t),\xi,b_2)$ admits a maximal solution $(\hat{y},\hat{z},\hat{k})$ such that $\hat{y}\in{\mathcal{S}}^\infty$. Moreover, for each $l\in{\mathcal{S}}^\infty$ such that $l_T\leq \xi$, Theorem \ref{th4.2}(i) also implies that the RBSDE$(G(t,y,z,\hat{Y}_t,Y_t),\xi,b_2\vee l)$ admits a solution $(y^l,z^l,k^l)$ such that $y^l\in{\mathcal{S}}^\infty$. Then, by Proposition \ref{pr3.3}(ii) with ${\mathcal{Y}}_D={\mathcal{S}}^\infty$, we deduce that the BSDE$(G(t,y,z,\hat{Y}_t,Y_t),\xi)$ admits a maximal solution $(\tilde{y},\tilde{z})$ such that $\tilde{y}\in{\mathcal{S}}^\infty$. The arguments above imply that for each $\theta\in(0,1)$, $\tilde{y}_t\geq\hat{y}^\theta_t\geq\frac{\hat{Y}_t-\theta Y_t}{1-\theta}$. Since $\tilde{y}_t(1-\theta)\geq \hat{Y}_t-\theta Y_t$, sending $\theta$ to $1$, we get $\hat{Y}_t\leq Y_t.$ Similarly, by considering the difference $\frac{Y_t-\theta\hat{Y}_t}{1-\theta}$ as in (\ref{4.8}), we can also get $Y_t\leq\hat{Y}_t.$ Thus, $(Y,Z)$ is a unique solution to the BSDE$(g,\xi)$ such that ${Y\in\mathcal{S}}^\infty$.

\textbf{Case (ii):} Let (\hyperref[4A1]{4A1}-(i)(iii)) hold and $D=(0,\infty)$. It can be checked that the backward ODE:
\begin{equation*}
\psi_{b_3}(t)=b_3+\int_t^Tu(s)\psi_{b_3}(s)\ln(\psi_{b_3}(s))ds,\ \  t\in[0,T]
\end{equation*}
admits a solution $\psi_{b_3}(t)=b_3^{\exp(\int_t^Tu(s)ds)}$. Since $0<b_3<1$, we have $0<\psi_{b_3}(t)\leq b_3$. By (\hyperref[4A1]{4A1}-(i)(iii)) and (\ref{4.7}), we have  $g(t,\varphi(t),0)\leq u(t)l(\varphi(t))$ and
\begin{equation*}
g(t,\psi_{b_3}(t),0)\geq -u(t)|\psi_{b_3}(t)\ln(|\psi_{b_3}(t)|)|=u(t)\psi_{b_3}(t)\ln(\psi_{b_3}(t)),\ \ \ dt\times dP\textrm{-}a.e.,
\end{equation*}
which together with the facts that $b_3\leq\xi\leq c_2$ and $\psi_{b_3}(t)\leq b_3\leq c_2\leq\varphi(t)$, implies that the BSDE$(g,\xi)$ is dominated by $\psi_{b_3}(t)$ and $\varphi(t)$. Then by Lemma \ref{l2.4}, the BSDE$(g,\xi)$ admits a minimal solution $(Y,Z)$ such that $Y_t\geq \psi_{b_3}(t)$ and $Y\in{\mathcal{S}}_D^\infty$.

Assume that (\hyperref[4A2]{4A2}) holds. Let $(\hat{Y},\hat{Z})$ be another solution to the BSDE$(g,\xi)$ such that $\hat{Y}_t\geq \psi_{b_3}(t)$ and ${\hat{Y}\in\mathcal{S}}_D^\infty$. Since $(Y,Z)$ is a minimal solution to the BSDE$(g,\xi)$ such that $Y_t\geq \psi_{b_3}(t)$ and ${Y\in\mathcal{S}}_D^\infty$, we have $Y_t\leq\hat{Y}_t$. Then by the proof of Theorem \ref{th4.2}(ii), we deduce that $Y_t\geq\hat{Y}_t$, which implies that $(Y,Z)$ is a unique solution to the BSDE$(g,\xi)$ such that $Y_t\geq \psi_{b_3}(t)$ and ${Y\in\mathcal{S}}_D^\infty$. From the arguments above, it follows that for each constant $\bar{b}$ such that $0<\bar{b}<b_3$, $(Y_t,Z_t)$ is actually a unique solution to the BSDE$(g,\xi)$ such that $Y_t\geq \psi_{\bar{b}}(t)$ and ${Y\in\mathcal{S}}_D^\infty$. Let $(\tilde{Y},\tilde{Z})$ be another solution to the BSDE$(g,\xi)$ such that $\tilde{Y}\in{\mathcal{S}}_D^\infty$. It follows that there exists a constant $0<\tilde{b}<b_3\wedge\tilde{Y}_t,$ $dt\times dP$-$a.e.,$ such that $\tilde{Y}_t\geq\psi_{\tilde{b}}(t)$. This means that $\tilde{Y}_t={Y}_t,$ and thus $(Y,Z)$ is a unique solution to the BSDE$(g,\xi)$ such that $Y_t\in{\mathcal{S}}_D^\infty$.
\end{proof}

\begin{remark}\label{r4.4} (i) The proofs of the uniqueness of solutions in Theorems \ref{th4.2} and \ref{th4.3} combine the $\theta$-difference technique from \cite{BH08} with some novel comparison arguments based on RBSDEs. This method differs from the related studies on quadratic RBSDEs in \cite{K2, BY, Li, Gu}. The use of comparison arguments is a key difference between our method and those based on the $\theta$-difference technique for quadratic BSDEs in \cite{BH08, Y, FH, FHT1, FHT2, LG}.

(ii) The existence results in Theorems \ref{th4.2} and \ref{th4.3} generalize the corresponding results in \cite{K2, Xu} and \cite{K, LS, BL, BE, F}, respectively. Some locally Lipschitz conditions were assumed in \cite{K2, BE, Li, Im} and (\hyperref[4A3]{4A3}) below to ensure the uniqueness of bounded solutions. However, the generator $g$ in our setting may not satisfy these conditions (see Example \ref{ex4.6}(ii) below). This constitutes an important feature of the one-sided conditions introduced in this paper. To still guarantee the uniqueness, we therefore impose (\hyperref[4A2]{4A2}) in Theorems \ref{th4.2} and \ref{th4.3}.
\end{remark}

Some special cases of the BSDEs and RBSDEs discussed in this paper have many important applications (see \cite{DE, SS, MZ, BT, Tian, LG, Zheng2, Zheng3}), which motivate our current investigation. We next present a motivating example from finance.

\begin{example}\label{ex4.5} Let $\eta\in \mathcal{S}^\infty_D$ be the value process of a portfolio and $f(x)\in C(D)$. The functional
\begin{equation*}
\rho_t(\eta_T):=u_f^{-1}(E[u_f(\eta_T)|{\cal{F}}_t]),\quad t\in[0,T],\footnote{For the definition and properties of the function $u_f(x)$, see Appendix \ref{appC}.}
\end{equation*}
can be considered as a dynamic risk measure (DRM) of $\eta_T$. In fact, if $f(x)$ is a constant, then by Example \ref{exC.2}(i), one can check that $\rho$ becomes the dynamic entropy risk measure (see Barrieu and El Karoui \cite{BK2}). If $f(x)$ is nonpositive on $D$, then $u_f(x)$ is a general utility function (see Lemma \ref{lmC.1}(ii) and Example \ref{exC.2}), and thus $\rho$ becomes a conditional certainty equivalent. By Zheng et al. \cite[Theorem 3.2]{Zheng2}, $\rho(\eta_T)$ is actually the first component of a unique solution to the BSDE$(f(y)|z|^2,\eta_T)$ such that $\rho(\eta_T)\in\mathcal{S}_D^\infty$.

Let $h(\omega,t,z): \Omega\times [0,T]\times{\mathbf{R}}^d\longmapsto{\mathbf{R}}$ be measurable with respect to ${\mathcal{P}}\otimes{\mathcal{B}}({\mathbf{R}})$ and sublinear in $z$ with $|h(\cdot,z)|\leq\kappa|z|$, $dt\times dP$-$a.e.$, where $\kappa>0$ is a constant. From the proof of Zheng \cite[Corollary 3.4]{Zheng4}, we derive that $Y_t:=\operatorname {ess\,sup}_{Q\in{\cal{Q}}^h}E_Q[u_f(\eta_T)|{\cal{F}}_t]$,\footnote{$\mathcal{Q}^h:=\{Q^q:$ $\frac{dQ^q}{dP}=\exp\{\int_0^Tq_s\cdot dB_s-\frac{1}{2}\int_0^T|q_s|^2ds\}$ with $q_t\in H^2_d$ and $dt\times dP$-$a.e.,$ $q_t\in\partial h(t,0)\}$, where $\partial h(t,0)$ is the subdifferential of $h(t,0)$.} is the first component of a unique solution to the BSDE$(h(t,z),u_f(\eta_T))$ such that $Y\in\mathcal{S}^\infty$. The functional
\begin{equation*}
\tilde{\rho}_t(\eta_T):=\operatorname {ess\,sup}_{Q\in{\cal{Q}}^h}u_f^{-1}(E_Q[u_f(\eta_T)|{\cal{F}}_t]),
\quad t\in[0,T],
\end{equation*}
is a robust DRM. By Lemma \ref{lmC.1}(iv), we have $\tilde{\rho}_t(\eta_T)=u_f^{-1}(Y_t)$. Then, applying It\^{o}'s formula to $u_f^{-1}(Y_t)$, and by Lemma \ref{lmC.1}(iii), we deduce that $\tilde{\rho}(\eta_T)$ is the first component of a solution to the BSDE$(H(t,y,z),\eta_T)$ such that $\tilde{\rho}(\eta_T)\in\mathcal{S}_D^\infty$, where $H(t,y,z):=\frac{h(t,u_f'(y)z)}{u_f'(y)}+f(y)|z|^2.$ In particular, if $D=(0,\infty)$, $f(y)=\frac{r-1}{2y}, r>1$, and $dt\times dP$-$a.e.$, $h(t,\cdot)\geq0$, then by Example \ref{exC.2}(ii), $\tilde{\rho}_t$ becomes the robust $r$-norm, which was given for unbounded variables by Laeven et al. \cite[Examples 42]{LG} via so-called Geometric BSDEs. In fact, the unbounded robust DRM $\tilde{\rho}_t$ (resp. robust $r$-norm) can also be constructed similarly, since we only require $u_f(\eta_T)\in L^p(\mathcal{F}_T)$ in the arguments above.

Furthermore, the functional
\begin{equation*}
\bar{\rho}_t(\eta_T):=\operatorname {ess\,sup}_{\tau\in{\cal{T}}_{t,T}}\operatorname {ess\,sup}_{Q\in{\cal{Q}}^h}u_f^{-1}(E_Q[u_f(\eta_\tau)|{\cal{F}}_t]),\quad t\in[0,T],\footnote{$\mathcal{T}_{t,T}$ consists of all the stopping times $\tau$ such that $t\leq\tau\leq T$.}
\end{equation*}
is a maximal robust DRM of $\eta$ over $[t,T]$. By \cite[Corollary 3.4]{Zheng4} again and the fact that $h(t,0)\equiv0$, $dt\times dP$-$a.e.$, we deduce that for each $\tau\in{\cal{T}}_{t,T},$ $Y_t^\tau:=\operatorname {ess\,sup}_{Q\in{\cal{Q}}^h}E_Q[u_f(\eta_\tau)|{\cal{F}}_t]$ is the first component of a unique solution to the BSDE$(h(t,z),u_f(\eta_\tau))$ such that $Y^\tau\in\mathcal{S}^\infty$ with $Y_s^\tau=u_f(\eta_\tau)$ for each $s\in[\tau,T]$. By Kobylanski et al. \cite[Proposition 3.1]{K2}, we get that $\hat{Y}_t:=\operatorname {ess\,sup}_{\tau\in{\cal{T}}_{t,T}}Y^\tau_t$ is the first component of a unique solution to the RBSDE$(h(t,z),u_f(\eta_T),u_f(\eta))$ such that $\hat{Y}\in\mathcal{S}^\infty$. By the arguments above and Lemma \ref{lmC.1}(iv), we have $\bar{\rho}_t(\eta_T)=u_f^{-1}(\hat{Y}_t)$. Then, applying It\^{o}'s formula to $u_f^{-1}(\hat{Y}_t)$ and by Lemma \ref{lmC.1}(iii), we get that $\bar{\rho}(\eta_T)$ is the first component of a solution to the RBSDE$(H(t,y,z),\eta_T,\eta)$ such that $\bar{\rho}(\eta_T)\in\mathcal{S}_D^\infty$.

Observe that $|H(t,y,z)|\leq\kappa|z|+|f(y)||z|^2$. Inspired by \cite{Peng, BK2, LG, Zheng4}, we can further develop general DRMs theory, using the well-posedness results of RBSDEs and BSDEs obtained in this subsection and Section \ref{sec5}.
\end{example}

To the best of our knowledge, the well-posedness results obtained in Theorems \ref{th4.2} and \ref{th4.3} have not been obtained in the literature. We show some examples.

\begin{example}\label{ex4.6} (i) Let $D=\bf{R}$, $\xi\in L^\infty({\mathcal{F}}_T)$, $L\in{\mathcal{S}}^\infty$, and
\begin{equation*}
\mathfrak{g}(t,y,z)=u(t)l(y)h(y)+\phi(t,y)k(y)+\psi(t,y)|z|^2,\quad (y,z)\in D\times\mathbf{R}^d,
\end{equation*}
where $u(t)\in C_+([0,T])$, $l(x)\in{\mathcal{L}}$, $h(y)$, $k(y)\in C(\mathbf{R})$, and $\phi(\omega,t,y),\psi(\omega,t,y): \Omega\times [0,T]\times{\mathbf{R}}\longmapsto{\mathbf{R}}$ are both measurable with respect to ${\mathcal{P}}\otimes{\mathcal{B}}({\mathbf{R}})$ and continuous on $[0,T]\times {\mathbf{R}}$. By Remark \ref{r2.1}, $\mathfrak{g}$ satisfies (\hyperref[2A1]{2A1}). If there exists a constant $c>0$ such that for each $y\geq c$, $h(y)\leq1, k(y)=0,$ then Theorem \ref{th4.2} implies that the RBSDE$(\mathfrak{g},\xi, L)$ admits a minimal solution $(Y,Z,K)$ such that $Y\in{\mathcal{S}}^\infty$; and moreover, if we also have for each $y\leq -c$, $h(y)\geq-1$, $k(y)=0,$ then Theorem \ref{th4.3} implies that the BSDE$(\mathfrak{g},\xi)$ admits a solution $(\hat{Y},\hat{Z})$ such that $\hat{Y}\in{\mathcal{S}}^\infty.$

(ii) Let $D=(0,\infty)$, $\xi\in L_D^\infty({\mathcal{F}}_T)$, $L\in{\mathcal{S}}_D^\infty$, and
\begin{equation*}
\mathfrak{g}(t,y,z)=b_1+\frac{b_2}{y^k}+b_3y\ln(y)+b_4\phi(t,y)+b_5|z|^l+\frac{b_6}{y^r}|z|^2+g_1(z),
\quad (y,z)\in D\times\mathbf{R}^d,
\end{equation*}
with
\begin{equation*}\label{4.10}
\phi(t,y)=\left\{
                 \begin{array}{ll}
                   (y-1)^2|\eta_t|+\exp(-|\eta_t|), & y<1;\\
                   \exp(-|\eta_t|), & y\geq 1,
                 \end{array}
               \right.\tag{4.10}
\end{equation*}
where $\eta\in\mathcal{C}$, and $k>0$, $1\leq l\leq2$, $0\leq r\leq1$ and $b_i\geq0$ $(1\leq i\leq6)$ are all constants, $g_1(z):{\bf{R}}^d\rightarrow{\bf{R}}$ is a bounded Lipschitz function with a bounded support and $g(0)\geq0$. By Remark \ref{r2.1}, $\mathfrak{g}$ satisfies (\hyperref[2A1]{2A1}) and (\hyperref[4A1]{4A1}-(i)(iii)). Since $\frac{1}{y^r}|z|^2$ is convex on $D\times\mathbf{R}^d$ (see Lemma \ref{lmD.1}), by Remark \ref{r4.1}(i)(iii), we can further get that $\mathfrak{g}$ satisfies (\hyperref[4A2]{4A2}). Hence, Theorem \ref{th4.2} implies that the RBSDE$(\mathfrak{g},\xi, L)$ admits a unique solution $(Y,Z,K)$ such that $Y\in{\mathcal{S}}_D^\infty$, and Theorem \ref{th4.3} implies that the BSDE$(\mathfrak{g},\xi)$ admits a unique solution $(\hat{Y},\hat{Z})$ such that $\hat{Y}\in{\mathcal{S}}_D^\infty.$ Since $\eta$ in (\ref{4.10}) belongs to $\mathcal{C}$, one can verify that $\mathfrak{g}$ does not satisfy the locally Lipschitz conditions assumed in \cite{K2, BE, Li, Im} or in (\hyperref[4A3]{4A3}) below. Those conditions are typically used to ensure the uniqueness of bounded solutions.
\end{example}
\subsection{Case of a general growth in $y$}\label{sec4.2}
Note that in this subsection, $D$ is an arbitrary open interval unless we specify it. We first consider the following two examples.
\begin{example}\label{ex4.7} (i) By El Karoui et al. \cite[Proposition 2.2]{EPQ} and Lemma \ref{lmB.1}, we get that the BSDE$(y,-\frac{1}{2})$ admits a unique solution $(\hat{Y},\hat{Z})$ such that $\hat{Y}\in{\mathcal{S}}^\infty$, and moreover, $\hat{Y}_t=-\frac{1}{2}e^{(T-t)}$. When $e^T<2$, by setting $Y_t:=-\ln(\hat{Y}_t+1)$ and $Z_t:=0$, we get that the BSDE$(-1+\exp(y)-\frac{1}{2}|z|^2,-\ln(\frac{1}{2}))$ admits a solution $(Y,Z)$ such that $Y\in{\mathcal{S}}^\infty$. But when $e^T\geq2$, if this BSDE admits a solution $(Y,Z)$ such that $Y\in{\mathcal{S}}^\infty$, then by applying It\^{o}'s formula to $\exp(-Y_t)-1$ and setting $\hat{y}_t:=\exp(-Y_t)-1$ and $\hat{z}_t:=-\exp(-Y_t)Z_t$, we get that $(\hat{y},\hat{z})$ is a unique solution to the BSDE$(y,-\frac{1}{2})$ such that $\hat{y}_0=-\frac{1}{2}e^T\leq-1$, which contradicts $\hat{y}_0=\exp(-Y_0)-1>-1$.

(ii) Let $D=(0,\infty)$ and $\beta>0$ be a constant. When $T<\frac{1}{\beta}$, the BSDE$(-\beta, 1)$ admits a unique solution $(\hat{Y},\hat{Z})$ such that $\hat{Y}\in{\mathcal{S}}_D^\infty$ with $\hat{Y}_t=1-\beta(T-t)$. But when $T\geq\frac{1}{\beta}$, this BSDE has no solution $(Y,Z)$ such that $Y\in{\mathcal{S}}_D^\infty$.
\end{example}
Example \ref{ex4.7} implies that the growth conditions on $y$ in (\hyperref[4A1]{4A1}-(i)(iii)) can be further extended when $T$ is small enough. Inspired by this, we will consider the existence of solutions under a general growth condition on $y$ for small time duration.
\begin{itemize}
  \item \textbf{(4A1')}\label{4A1'} There exist a function $l(x,y)\in C_+(\mathbf{R}\times D)$ and an interval $[b,c]\subset D$, such that for each $(t,y)\in[0,T]\times D,$
\begin{itemize}
  \item\textbf{(i)}  \ \ \ \ \ \ \ \ \ \ \ \ \ \ \ \ \ \ \ \ \ \ \ \ \ \ \ \ \ \ \ \ \  $1_{\{y\geq c\}}{g}(t,y,0)\leq l(t,y);$
  \item\textbf{(ii)} \ \ \ \ \ \ \ \ \ \ \ \ \ \ \ \ \ \ \ \ \ \ \ \ \ \ \ \ \ \ \ $1_{\{y\leq b\}}{g}(t,y,0)\geq -l(t,y).$
\end{itemize}
\end{itemize}

By Peano existence theorem, we deduce that for each $a\in D$ and $l(x,y)\in C_+(\mathbf{R}\times D)$, there exists a constant $\lambda>0$, such that whenever $T\leq\lambda$, the backward ODE:
\begin{equation*}
\varphi(t)=a+\int_t^Tl(s,\varphi(s))ds, \ \ t\in[0,T]
\end{equation*}
admits a solution $\varphi(t)\in C(D).$ Then, by similar arguments as in the proofs of Theorem \ref{th4.2} and \ref{th4.3}, we obtain the following existence results:
\begin{proposition}\label{pr4.8}
Let (\hyperref[4A1']{4A1'}-(i)) hold, $\xi\in L_D^\infty({\mathcal{F}}_T)$, and $L\in{\mathcal{S}}_D^\infty$. Then there exists a constant $\lambda>0$, such that whenever $T\leq\lambda$, the RBSDE$(g,\xi,L)$ admits a minimal solution $(Y,Z,K)$ such that ${Y\in\mathcal{S}}_D^\infty$.
\end{proposition}

\begin{proposition}\label{pr4.9}
Let (\hyperref[4A1']{4A1'}-(i)(ii)) hold and $\xi\in L_D^\infty({\mathcal{F}}_T)$. Then there exists a constant $\lambda>0$, such that whenever $T\leq\lambda$, the BSDE$(g,\xi)$ admits at least one solution $(Y,Z)$ such that ${Y\in\mathcal{S}}_D^\infty$.
\end{proposition}

\begin{remark}\label{r4.10} By Proposition \ref{pr4.9} and Lemma \ref{lmB.1}, it can be checked that the ‘‘characteristic BSDE" in \cite[Equation (3.8)]{MZ}) with the terminal variable $\xi\in L_D^\infty({\mathcal{F}}_T)$ admits at least one solution $(Y,Z)\in{\mathcal{S}}_D^\infty\times {\mathcal{H}}^{BMO}_d$, when the process $\sigma_3$ therein has a lower bound $M>0$, $D=(\frac{1}{M}, \infty)$ and $T$ is small enough.
\end{remark}

To obtain a uniqueness result for the BSDE with a general growth in $y$, we introduce the following locally Lipschitz condition:
\begin{itemize}
  \item \textbf{(4A3)}\label{4A3} There exist two nonnegative processes $r^1, r^2\in{\mathcal{H}}^{BMO}_1$, a constant $r\in(0,1)$ and two continuous functions $h_1(x,y),h_2(x,y):D\times D\mapsto[0,\infty)$, such that $dt\times dP\textrm{-}a.e.,$ for each $y_1, y_2\in D$ and $z_1,z_2\in\mathbf{R}^d$, $$|{g}(t,y_1,z_1)-{g}(t,y_2,z_2)|\leq \hat{h}_1|y_1-y_2|+\hat{h}_2|z_1-z_2|,$$
      with $\hat{h}_1:=h_1(|y_1|,|y_2|)(r_t^1+|z_1|^{2r}+|z_2|^{2r})$ and $\hat{h}_2:=h_2(|y_1|,|y_2|)(r_t^2+|z_1|+|z_2|).$
\end{itemize}

The following result can be seen as an extension of the uniqueness results in Briand and Elie \cite[Corollary 2.2]{BE} and Imkeller et al. \cite[Theorem 3.4]{Im}.

\begin{proposition}\label{pr4.11} Let (\hyperref[4A3]{4A3}) hold and $\sqrt{|g(\cdot,\beta,0)|}\in {\mathcal{H}}^{BMO}_1$ for some constant $\beta\in D$. Then the BSDE$(g,\xi)$ admits at most one solution $(Y,Z)$ such that $Y\in{\mathcal{S}}_D^{\infty}$.
\end{proposition}
\begin{proof} Let the BSDE$(g,\xi)$ admit a solution $(Y,Z)$ such that $Y\in{\mathcal{S}}_D^{\infty}$. When (\hyperref[4A3]{4A3}) holds, we have $dt\times dP$-$a.e.,$
\begin{align}
|g(t,Y_t,Z_t)|&\leq|g(t,\beta,0)|+|g(t,Y_t,0)-g(t,\beta,0)|+|g(t,Y_t,Z_t)-g(t,Y_t,0)|\tag{4.11}\label{4.11}\\
&\leq|g(t,\beta,0)|+h_1(|Y_t|,|\beta|)r_t^1|Y_t-\beta|+h_2(|Y_t|,|Y_t|)(r_t^2+|Z_t|)|Z_t|\notag\\
&\leq|g(t,\beta,0)|+h_1(|Y_t|,|\beta|)r_t^1|\beta|+h_2(|Y_t|,|Y_t|)|r_t^2|^2\notag\\
&\ \ \ +h_1(|Y_t|,|\beta|)r_t^1|Y_t|+2h_2(|Y_t|,|Y_t|)|Z_t|^2.\notag
\end{align}
Since $\sqrt{|g(\cdot,\beta,0)|}\in {\mathcal{H}}^{BMO}_1$, by (\ref{4.11}), there exist two nonnegative processes $\sqrt{\eta}\in{\mathcal{H}}^{BMO}_1$ and $C\in{\mathcal{S}}^\infty$, such that $g(t,Y_t,Z_t)\leq\eta_t+C_t|Z_t|^2$, $dt\times dP$-$a.e.$ Then by Lemma \ref{lmB.1}, we have $Z\in{\mathcal{H}}^{BMO}_d.$

Let $(\hat{Y},\hat{Z})\in {\mathcal{S}}_D^{\infty}\times {\mathcal{H}}^{BMO}_d$ be another solution to the BSDE$(g,\xi)$. Then by a linearization argument, we have
\begin{equation}\label{4.12}
  \hat{Y}_t-Y_t=0+\int_t^T(b_s(\hat{Y}_s-Y_s)+c_s\cdot (\hat{Z}_s-Z_s))ds-\int_t^T(\hat{Z}_s-Z_s)\cdot dB_s,\ \ t\in[0,T], \tag{4.12}
\end{equation}
where
\begin{equation*}
b_s=\frac{g(s,\hat{Y}_s,\hat{Z}_s)-g(s,Y_s,\hat{Z}_s)}{\hat{Y}_s-Y_s}1_{\{|\hat{Y}_s-Y_s|>0\}} \ \ \textrm{and}\ \
c_s=\frac{g(s,Y_s,\hat{Z}_s)-g(s,Y_s,Z_s)}{|\hat{Z}_s-Z_s|^2}(\hat{Z}_s-Z_s)1_{\{|\hat{Z}_s-Z_s|>0\}}.
\end{equation*}
Set $\tilde{G}(t,y,z):=b_ty+c_tz.$ Since $|b_t|\leq h_1(|\hat{Y}_t|,|Y_t|)(r_t^1+|\hat{Z}_t|^{2r}+|Z_t|^{2r})$ and $|c_t|\leq h_2(|\hat{Y}_t|,|Y_t|)(r_t^2+|\hat{Z}_t|+|Z_t|)$, by Briand and Confortola \cite[Theorem 10]{BC}, the BSDE$(\tilde{G},0)$ admits a unique solution $(0,0)\in{\mathcal{S}}^\infty\times{\mathcal{H}}_d^{BMO}$. This with (\ref{4.12}) gives $(\hat{Y}_t-Y_t,\hat{Z}_t-Z_t)=(0,0)$, $dt\times dP\textrm{-}a.e.$ The proof is complete.
\end{proof}
\section{The existence and uniqueness of unbounded solutions}\label{sec5}
\subsection{Quadratic reflected BSDEs}
We introduce the following assumptions for $g$:
\begin{itemize}
  \item \textbf{(5A1)}\label{5A1} There exist a function $f(x)\in C_+(D)$, three nonnegative processes $\delta\in H_1^1, \gamma,\vartheta\in{\mathcal{C}}$, and four constants $\kappa\geq0,$ $\nu\geq0,$ $c\geq b>0,$ such that $dt\times dP\textrm{-}a.e.,$ for each $(y,z)\in D\times{\mathbf{R}}^{\mathit{d}},$
      \begin{itemize}
  \item\textbf{(i)}  \ \ \ \ \ \ \ \ \ \ \ \ \ \ \ \ \ \ \ \ \ $1_{\{y\geq c\}}{g}(t,y,z)\leq\delta_t+\gamma_t |y|+\kappa|z|+f(y)|z|^2;$
  \item\textbf{(ii)} \ \ \ \ \ \ \ \ \ \ \ \ \ \ \ \ \ \ $1_{\{y\leq -c\}}{g}(t,y,z)\geq-\delta_t-\gamma_t|y|-\kappa|z|-f(-y)|z|^2;$
  \item\textbf{(iii)} \ \ \ \ \ \ \ \ \ \ \ \ \ \ \ \ \ \ \ \ \ \ \ \ \ $1_{\{y\leq b\}}{g}(t,y,z)\geq-\vartheta_t|y|-\kappa|z|-\frac{\nu}{y}|z|^2.$
\end{itemize}
  \item \textbf{(5A2)}\label{5A2} $g$ satisfies the $\theta$-domination condition for ${\mathcal{C}}_D$ and $H^2_d$ with $G$, and there exist two nonnegative processes $\tilde{\mu}\in H_1^1, \tilde{\gamma}\in{\mathcal{C}}$ and a constant $\tilde{\kappa}\geq0$, such that $dt\times dP\textrm{-}a.e.,$ for any $(y,z)\in D\times{\mathbf{R}}^{\mathit{d}}$ and $y^1,y^2\in{\mathcal{C}}_D$,
\begin{equation*}
G(t,y,z,y_t^1,y_t^2)\leq\tilde{\mu}_t(1+|y_t^1|+|y_t^2|)+\tilde{\gamma}_t |y|+\tilde{\kappa}|z|+{g}(t,y,z).
\end{equation*}
  \item \textbf{(5A2')}\label{5A2'} $dt\times dP\textrm{-}a.e.,$ for each $y\in D$, $g(t,y,\cdot)$ is convex in $z$, and there exists a nonnegative process ${\mu}\in{\mathcal{C}}$, such that for each $y_1,y_2\in D$ and $z\in\mathbf{R}^d$,
      \begin{equation*}
      |{g}(t,y_1,z)-{g}(t,y_2,z)|\leq \mu_t|y_1-y_2|.
      \end{equation*}
\end{itemize}

\begin{remark}\label{r5.1}
\begin{itemize}
\item (i) From Remark \ref{r4.1}(i) and (\ref{4.2}), it follows that (\hyperref[4A2']{4A2'}) and (\hyperref[5A2']{5A2'}) are both the special cases of (\hyperref[5A2]{5A2}). (\hyperref[5A1]{5A1}-(i)) combined with (\hyperref[5A2]{5A2}) is similar to \cite[(H2')]{FH}. However, a key distinction lies in the treatment of singular generators: \cite[(H2')]{FH} can not contain the singular terms such as $\frac{|z|^2}{y}$ and $\frac{1}{y^2}$, etc.
\item (ii) Let $D=\mathbf{R}$. Even when $\gamma\equiv\alpha$ and $f(y)\equiv\beta,$ where $\alpha$ and $\beta$ are two nonnegative constants,  (\hyperref[5A1]{5A1}-(i)(ii)) are still more general than the one-sided conditions for quadratic BSDEs in \cite{FHT1, FHT2}:
\begin{equation*}\label{eq5.1}
\textrm{sgn}(y){g}(t,y,z)\leq\delta_t+\alpha|y|+\beta|z|^2.\tag{5.1}
\end{equation*}
Indeed, in this case, if $g(t,y,z)=\phi(t,y)|z|^2$, where $\phi(t,y)$ is the function defined in (\ref{4.10}), then $g$ satisfies (\hyperref[5A1]{5A1}-(i)(ii)), but does not satisfy (\ref{eq5.1}).
\end{itemize}
\end{remark}

For convenience, for $\delta, \gamma\in H_1^1$ and $t\in[0,T]$, we set
\begin{equation*}
\Lambda^{\delta,\gamma}_t(\zeta):=e^{\int_0^t\gamma_sds}\left(|\zeta|+\int_0^t\delta_sds\right),\quad \forall \zeta\in L({\mathcal{F}}_t).
\end{equation*}
and denote $(\Lambda^{\delta,\gamma}_t(\eta_t))_{t\in[0,T]}$ by $\Lambda^{\delta,\gamma}(\eta)$ for $\eta\in H_1^1$. Then, we have the following existence results for unbounded solutions to RBSDEs. A related result was obtained by \cite[Theorem 3.1]{BY}, where $|g|\leq\beta(1+|y|+|z|^2)$, where $\beta>0$ is a constant.

\begin{proposition}\label{pr5.2} Let (\hyperref[5A1]{5A1}-(i)) hold with $f(y)$ nondecreasing. Let $u_f(\Lambda^{\delta,\gamma}_T(\xi\vee c))\in L^p({\mathcal{F}}_T)$ and $u_f(\Lambda^{\delta,\gamma}(L))\in{\mathcal{S}}^p$. Then the RBSDE$(g,\xi,L)$ admits:

(i) A minimal solution $(\underline{Y},\underline{Z},\underline{K})$ such that $\underline{Y}\in\mathcal{C}_D$, in particular $u_f(\Lambda^{\delta,\gamma}(\underline{Y}))\in{\mathcal{S}}^p$;

(ii) A maximal solution $(\overline{Y},\overline{Z},\overline{K})$ such that $u_f(\Lambda^{\delta,\gamma}(\overline{Y}))\in{\mathcal{S}}^p$.
\end{proposition}
\begin{proof} We divide this proof into four steps.

\textbf{Step 1.} For $(\omega,t,y,z)\in \Omega\times[0,T]\times D\times{\mathbf{R}}^{\mathit{d}}$, we set
\begin{align*}\label{5.2}
g_1(\omega,t,y,z)&:=\kappa|z|+f(y)|z|^2;\\
g_2(\omega,t,y,z)&:=\kappa|z|+f((y\vee\Lambda_t^{\delta,\gamma}(c))-e^{\int_0^t\gamma_sds}\int_0^t\delta_sds)|z|^2;\tag{5.2}\\
g_3(\omega,t,y,z)&:=\gamma_t e^{\int_0^t\gamma_sds}\int_0^t\delta_sds+e^{\int_0^t\gamma_sds}\delta_t+\kappa|z|+f(y)|z|^2;\\ g_4(\omega,t,y,z)&:=e^{\int_0^t\gamma_sds}\delta_t+\kappa|z|+e^{-\int_0^t\gamma_sds}f(e^{-\int_0^t\gamma_sds}y)|z|^2;
\\g_5(\omega,t,y,z)&:=\delta_t+\gamma_ty+\kappa|z|+f(y)|z|^2.
\end{align*}
By Remark \ref{r2.1}, we get that all the generators $g_1$-$g_5$ satisfy (\hyperref[2A1]{2A1}).

\textbf{Step 2.} Since for each $t\in[0,T]$,
\begin{equation*}
u_f(c)\leq u_f(\Lambda^{\delta,\gamma}_t(L_t\vee c))\leq u_f(\Lambda^{\delta,\gamma}_t(L_t))+u_f(\Lambda^{\delta,\gamma}_T(\xi\vee c)),
\end{equation*}
we have $u_f(\Lambda^{\delta,\gamma}(L\vee c))\in{\mathcal{S}}^p$. By \cite[Corollary 3.5]{Zheng3}, the RBSDE$(g_1,\Lambda^{\delta,\gamma}_T(\xi\vee c),\Lambda^{\delta,\gamma}(L\vee c))$ admits a unique solution $(\bar{y}^1,\bar{z}^1,\bar{k}^1)$ such that $u_{f}(\bar{y}^1)\in{\mathcal{S}}^p$. Similarly, for each $l\in\mathcal{C}_D$ such that $u_{f}(l)\in{\mathcal{S}}^p$ and $l_T\leq\Lambda^{\delta,\gamma}_T(\xi\vee c)$, the RBSDE$(g_1,\Lambda^{\delta,\gamma}_T(\xi\vee c),\Lambda^{\delta,\gamma}(L\vee c)\vee l)$ admits a unique solution $(y^l,z^l,k^l)$ such that $u_{f}(y^l)\in{\mathcal{S}}^p$. Then by Proposition \ref{pr3.3}(i) with ${\mathcal{Y}}_D=\{\eta\in{\mathcal{C}}_D:u_{f}(\eta)\in{\mathcal{S}}^p\}$, the RBSDE$(g_2$, $\Lambda^{\delta,\gamma}_T(\xi\vee c),\Lambda^{\delta,\gamma}(L\vee c))$ admits a maximal solution $(\bar{y}^2,\bar{z}^2,\bar{k}^2)$ such that $u_{f}(\bar{y}^2)\in{\mathcal{S}}^p.$ Set
 \begin{equation*}\label{5.3}
\bar{y}^3_t:=\bar{y}^2_t-e^{\int_0^t\gamma_sds}\int_0^t\delta_sds,\ \bar{z}^3_t:=\bar{z}^2_t,\ \bar{k}^3_t:=\bar{k}^2_t,\ \ t\in[0,T].\tag{5.3}
\end{equation*}
For each $t\in[0,T]$, $\bar{y}^2_t\geq\Lambda^{\delta,\gamma}_t(c)$, by (\ref{5.2}), we have
\begin{equation*}
g_2(t,\bar{y}^2_t,\bar{z}^2_t)=\kappa|\bar{z}^2_t|+f(\bar{y}^2_t-e^{\int_0^t\gamma_sds}\int_0^t\delta_sds)|\bar{z}^2_t|^2.
\end{equation*}
Thus, by (\ref{5.3}) and It\^{o}'s formula, we can deduce that $(\bar{y}^3,\bar{z}^3,\bar{k}^3)$ is actually a maximal solution to the RBSDE$(g_3,e^{\int_0^T\gamma_sds}(\xi\vee c),e^{\int_0^\cdot\gamma_sds}(L\vee c))$ such that $u_{f}(\bar{y}^3+e^{\int_0^\cdot\gamma_sds}\int_0^\cdot\delta_sds)\in{\mathcal{S}}^p.$

\textbf{Step 3.} By the conclusion of Step 2, for each $l\in\mathcal{C}_D$ such that $u_{f}(l+e^{\int_0^\cdot\gamma_sds}\int_0^\cdot\delta_sds)\in{\mathcal{S}}^p$ and $l_T\leq e^{\int_0^T\gamma_sds}(\xi\vee c)$, we deduce that the RBSDE$(g_3, e^{\int_0^T\gamma_sds}(\xi\vee c), e^{\int_0^\cdot\gamma_sds}(L\vee c)\vee l)$ admits a solution $(\tilde{y}^l,\tilde{z}^l,\tilde{k}^l)$ such that $u_{f}(\tilde{y}^l+e^{\int_0^\cdot\gamma_sds}\int_0^\cdot\delta_sds)\in{\mathcal{S}}^p$. Then, by Proposition \ref{pr3.3}(i) with ${\mathcal{Y}}_D=\{\eta\in{\mathcal{C}}_D: u_{f}(\eta+e^{\int_0^\cdot\gamma_sds}\int_0^\cdot\delta_sds)\in{\mathcal{S}}^p\}$, the RBSDE$(g_4,e^{\int_0^T\gamma_sds}(\xi\vee c),e^{\int_0^\cdot\gamma_sds}(L\vee c))$ admits a maximal solution $(\bar{y}^4,\bar{z}^4,\bar{k}^4)$ such that $u_{f}(\bar{y}^4+e^{\int_0^\cdot\gamma_sds}\int_0^\cdot\delta_sds)\in{\mathcal{S}}^p.$ Set
\begin{equation*}
\bar{y}^5_t:=e^{-\int_0^t\gamma_sds}\bar{y}^4_t,\ \bar{z}^5_t:=e^{-\int_0^t\gamma_sds}\bar{z}^4_t,\ \bar{k}^5_t:=\int_0^te^{-\int_0^s\gamma_rdr}d\bar{k}^4_s,\ \ t\in[0,T].
\end{equation*}
By It\^{o}'s formula, we deduce that $(\bar{y}^5,\bar{z}^5,\bar{k}^5)$ is a maximal solution to the RBSDE$(g_5,\xi\vee c,L\vee c)$ such that $u_{f}(\Lambda^{\delta,\gamma}(\bar{y}^5))\in{\mathcal{S}}^p$.

\textbf{Step 4.} By (\hyperref[5A1]{5A1}-(i)), we get that $dt\times dP$-$a.s.$,
$g(t,\bar{y}^5_t,\bar{z}^5_t)\leq g_5(t,\bar{y}^5_t,\bar{z}^5_t).$
It follows that the RBSDE$(g,\xi,L)$ is dominated by $\bar{y}^5$. Then, by Lemma \ref{l2.3}, the RBSDE$(g,\xi,L)$ admits a minimal solution $(\underline{Y},\underline{Z},\underline{K})$ such that $\underline{Y}\in{\mathcal{C}}_D$. Since $L_t\leq\underline{Y}_t\leq\bar{y}^5_t$, we have
\begin{equation*}
u_0(0)\leq u_0(\Lambda^{\delta,\gamma}_t(\underline{Y}_t))\leq u_f(\Lambda^{\delta,\gamma}_t(\underline{Y}_t))\leq u_f(\Lambda^{\delta,\gamma}_t(L_t))\vee u_f(\Lambda^{\delta,\gamma}_t(\bar{y}^5)),
\end{equation*}
which implies $u_f(\Lambda^{\delta,\gamma}(\underline{Y}))\in{\mathcal{S}}^p$. We obtain (i).

By the conclusion of Step 3, for each $l\in\mathcal{C}_D$ such that $u_{f}(\Lambda^{\delta,\gamma}(l))\in{\mathcal{S}}^p$ and $l_T\leq\xi\vee c$, we deduce that the RBSDE$(g_5, \xi\vee c, L\vee c\vee l)$ admits a solution $(\hat{y}^l,\hat{z}^l,\hat{k}^l)$ such that $u_{f}(\Lambda^{\delta,\gamma}(\hat{y}^l))\in{\mathcal{S}}^p$. Then, by Proposition \ref{pr3.3}(i) with ${\mathcal{Y}}_D=\{\eta\in{\mathcal{C}}_D:u_{f}(\Lambda^{\delta,\gamma}(\eta))\in{\mathcal{S}}^p\}$, we get that the RBSDE$(g,\xi,L)$ admits a maximal solution $(\overline{Y},\overline{Z},\overline{K})$ such that  $u_{f}(\Lambda^{\delta,\gamma}(\overline{Y}))\in{\mathcal{S}}^p$. We obtain (ii).
\end{proof}

\begin{proposition}\label{pr5.3} Let (\hyperref[5A1]{5A1}-(i)) hold with $f(y)$ nondecreasing. Let $u_{f}(q\Lambda^{\delta, \gamma}_T(\xi\vee c))\in L^p({\mathcal{F}}_T)$ and $u_{f}(q\Lambda^{\delta, \gamma}(L))\in{\mathcal{S}}^p$ for each $q\geq1$. If one of the following two conditions holds:

(i) $\int_0^T\gamma_tdt\in L^\infty(\mathcal{F}_T)$ and (\hyperref[5A2]{5A2}) holds with $\int_0^T(\tilde{\mu}_t+\tilde{\gamma}_t)dt\in L^\infty(\mathcal{F}_T)$;

(ii) (\hyperref[4A2']{4A2'}) holds,
\\
then the RBSDE$(g,\xi,L)$ admits a unique solution $(Y,Z,K)$ such that $u_{f}(q\Lambda^{\delta, \gamma}(Y))\in{\mathcal{S}}^p$ for each $q\geq1$, in particular,
\begin{equation*}\label{5.4}
\left\{
  \begin{array}{ll}
    \exp(\Lambda^{\delta,\gamma}(Y))\in\bigcap_{q\geq1}{\mathcal{S}}^q,\ & when\ f\equiv\beta>0;\\
    \Lambda^{\delta,\gamma}(Y)\in{\mathcal{S}}^p, & when\ f\equiv0.
  \end{array}
\right.\tag{5.4}
\end{equation*}
\end{proposition}
\begin{proof}
For each $q\geq1$, we set $f^q(x):=qf(qx), x\in D.$ From the definition of $u_f$ (see Appendix \ref{appC}), it can be checked that for each $q\geq1$, there exist two constants $b_1$ and $b_2$ such that for each $x\in D,$
\begin{align*}\label{5.5}
u_{{f}}(qx)&=\int_{q}^{qx}\exp\left(2\int_{1}^{y}f(z)dz\right)dy
+\int_{1}^{q}\exp\left(2\int_{1}^{y}f(z)dz\right)dy\\
&=q\int_{1}^{x}\exp\left(2\int_{1}^{qy}f(z)dz\right)dy
+b_2\\
&=q\exp\left(2\int_{1}^{q}{f}(z)dz\right)\int_{1}^{x}\exp\left(2\int_{1}^{y}qf(qz)dz\right)dy
+b_2\\&
=b_1u_{{f}^q}(x)+b_2
\tag{5.5}
\end{align*}
with $b_1=q\exp(2\int_{1}^{q}{f}(z)dz)$ and $b_2=\int_{1}^{q}\exp(2\int_{1}^{y}{f}(z)dz)dy.$ This implies that for each $q\geq1$, $u_{f^q}(\Lambda^{\delta, \gamma}_T(\xi\vee c))\in L^p({\mathcal{F}}_T)$ and $u_{f^q}(\Lambda^{\delta, \gamma}(L))\in{\mathcal{S}}^p$. Then, by Proposition \ref{pr5.2}(i) and (\ref{5.5}), we deduce that the RBSDE$(g,\xi,L)$ admits a minimal solution $(Y,Z,K)$ such that for each $q\geq1$, $u_{f}(q\Lambda^{\delta,\gamma}(Y))\in{\mathcal{S}}^p$.

We next prove the uniqueness. Let $(\hat{Y},\hat{Z},\hat{K})$ be another solution to the RBSDE$(g,\xi,L)$ such that for each $q\geq1$, $u_{f}(q\Lambda^{\delta, \gamma}_t(\hat{Y}_t))\in{\mathcal{S}}^p$. For $\theta\in(0,1)$, we have
\begin{eqnarray*}
\frac{\hat{Y}_t-\theta {Y}_t}{1-\theta}&=&\xi+\int_t^T\frac{\hat{g}_s}{1-\theta}ds+\frac{(\hat{K}_T-\theta {K}_T)-(\hat{K}_t-\theta {K}_t)}{1-\theta}-\int_t^T\frac{\hat{Z}_s-\theta {Z}_s}{1-\theta}\cdot dB_s
\end{eqnarray*}
with $\frac{\hat{g}_s}{1-\theta}:=\frac{1}{{1-\theta}}(g(s,\hat{Y}_s,\hat{Z}_s)-\theta g(s,{Y}_s,{Z}_s)).$

\textbf{Case (i):} By (\hyperref[5A2]{5A2}) and (\hyperref[5A1]{5A1}-(i)), we have $dt\times dP\textrm{-}a.e.,$ $\frac{\hat{g}_t}{1-\theta}\leq G(t,\frac{\hat{Y}_t-\theta {Y}_t}{1-\theta},\frac{\hat{Z}_t-\theta {Z}_t}{1-\theta},\hat{Y}_t,Y_t)$ and for each $(y,z)\in D\times{\mathbf{R}}^{\mathit{d}},$
\begin{equation*}\label{5.6}
1_{\{y\geq c\}}G(t,y,z,\hat{Y}_t,Y_t)\leq\bar{\delta}_t+\bar{\gamma}_t|y|
+\bar{\kappa}|z|+f(y)|z|^2,\tag{5.6}
\end{equation*}
where $\bar{\delta}_t:=\tilde{\mu}_t(1+|\hat{Y}_t|+|Y_t|)+\delta_t$, $\bar{\gamma}_t:=\tilde{\gamma}_t+\gamma_t$ and $\bar{\kappa}:=\tilde{\kappa}+\kappa.$ Since $f$ is nonnegative, we get that $u_f$ is convex (see Lemma \ref{lmC.1}(ii)). Then, by setting $M:=4\|e^{\int_0^T\bar{\gamma}_tdt}\|_\infty(\|\int_0^T\tilde{\mu}_tdt\|_\infty+1)$, we have
\begin{align*}
u_0(0)&\leq u_{f}(\Lambda^{\bar{\delta}, \bar{\gamma}}_T(\xi\vee c))\\&\leq u_{f}\left(e^{\int_0^T\bar{\gamma}_tdt}\Lambda^{\delta, \gamma}_T(\xi\vee c)+e^{\int_0^T\bar{\gamma}_tdt}\int_0^T\tilde{\mu}_tdt+e^{\int_0^T\bar{\gamma}_tdt}\int_0^T\tilde{\mu}_t|\hat{Y}_t|dt
+e^{\int_0^T\bar{\gamma}_tdt}\int_0^T\tilde{\mu}_t|{Y}_t|dt\right)\\&\leq \frac{1}{4}u_{f}(M\Lambda^{\delta, {\gamma}}_T(\xi\vee c))+\frac{1}{4}u_{f}(M)+\frac{1}{4}\sup_{t\in[0,T]}u_{f}\left(M|\hat{Y}_t|\right)
+\frac{1}{4}\sup_{t\in[0,T]}u_{f}\left(M|Y_t|\right)
\\&\leq\frac{1}{4}u_{f}(M\Lambda^{\delta, {\gamma}}_T(\xi\vee c))+\frac{1}{4}u_{f}(M)+\frac{1}{4}\sup_{t\in[0,T]}u_{f}\left(M\Lambda^{\delta, {\gamma}}_t(\hat{Y}_t)\right)+\frac{1}{4}\sup_{t\in[0,T]}u_{f}\left(M\Lambda^{\delta, {\gamma}}_t(Y_t)\right)
,\tag{5.7}\label{5.7}
\end{align*}
which implies $u_{f}(\Lambda^{\bar{\delta}, \bar{\gamma}}_T(\xi\vee c))\in L^p({\mathcal{F}}_T)$. Similarly, we can also get that $u_{f}(q\Lambda^{\bar{\delta}, \bar{\gamma}}(\hat{Y})), u_{f}(q\Lambda^{\bar{\delta}, \bar{\gamma}}(Y))\in{\mathcal{S}}^p$ for each $q\geq1$. This with the inequality that for each $\theta\in(0,1)$,
\begin{align*}\label{5.8}
u_0(0)\leq u_f\left(\Lambda^{\bar{\delta}, \bar{\gamma}}_t\left(\frac{\hat{Y}_t-\theta Y_t}{1-\theta}\right)\right)
&\leq\frac{1}{2}u_f\left(2\Lambda^{\bar{\delta}, \bar{\gamma}}_t\left(\frac{1}{1-\theta}|\hat{Y}_t|\right)\right)
+\frac{1}{2}u_f\left(2\Lambda^{\bar{\delta}, \bar{\gamma}}_t\left(\frac{\theta}{1-\theta}|Y_t|\right)\right)\\
&\leq\frac{1}{2}u_f\left(\frac{2}{1-\theta}\Lambda^{\bar{\delta}, \bar{\gamma}}_t(\hat{Y}_t)\right)
+\frac{1}{2}u_f\left(\frac{2}{1-\theta}\Lambda^{\bar{\delta}, \bar{\gamma}}_t(Y_t)\right),\tag{5.8}
\end{align*}
implies $u_f(\Lambda^{\bar{\delta}, \bar{\gamma}}(\frac{\hat{Y}-\theta Y}{1-\theta}))\in{\mathcal{S}}^p$. Hence, by (\ref{5.6})-(\ref{5.8}) and Proposition \ref{pr5.2}, we get that for each $\theta\in(0,1)$, the RBSDE$(G(t,y,z,\hat{Y}_t,Y_t),\xi,\frac{\hat{Y}-\theta Y}{1-\theta})$ admits a solution $(y^\theta,z^\theta,k^\theta)$ such that $u_f(\Lambda^{\bar{\delta}, \bar{\gamma}}(y^\theta))\in{\mathcal{S}}^p$, and that the RBSDE$(G(t,y,z,\hat{Y}_t,Y_t),\xi,\hat{Y}_t)$ admits a maximal solution $(\overline{Y},\overline{Z},\overline{K})$ such that $u_f(\Lambda^{\bar{\delta}, \bar{\gamma}}(\overline{Y}))\in{\mathcal{S}}^p$.

Since $Y_t\leq \hat{Y}_t$, by the Step 3 in the proof of Theorem \ref{th4.2}, we deduce that for each $\theta\in(0,1)$, $(y^\theta_t,z^\theta_t,k^\theta_t)$ is a solution to the RBSDE$(G(t,y,z,\hat{Y}_t,Y_t),\xi,\hat{Y}_t)$. This implies that for each $\theta\in(0,1)$, $\overline{Y}_t\geq\frac{\hat{Y}_t-\theta Y_t}{1-\theta}.$ By sending $\theta$ to $1$, we get $ \hat{Y}_t\leq Y_t,$ which implies $(\hat{Y}_t,\hat{Z}_t,\hat{K}_t)=(Y_t,Z_t,K_t).$

\textbf{Case (ii):} By (\hyperref[4A2']{4A2'}) and Remark \ref{r4.1}(i), we have $dt\times dP\textrm{-}a.e.,$ $\frac{\hat{g}_s}{1-\theta}\leq g(t,\frac{\hat{Y}_t-\theta{Y}_t}{1-\theta},\frac{\hat{Z}_t-\theta{Z}_t}{1-\theta})$.
Similar to (\ref{5.8}), we get that for each $\theta\in(0,1)$, $u_f(\Lambda^{\delta, \gamma}(\frac{\hat{Y}-\theta Y}{1-\theta}))\in{\mathcal{S}}^p$. Then, by a similar argument as in Case (i) (replacing $G$ therein with $g$), we can obtain $(\hat{Y}_t,\hat{Z}_t,\hat{K}_t)=(Y_t,Z_t,K_t).$

By Example \ref{exC.2}(i), we further obtain (\ref{5.4}). The proof is complete.
\end{proof}

\begin{remark}\label{r5.4} (i) Assume that $1_{\{y\geq c\}}{g}(t,y,z)\leq\delta_t+\gamma_t |y|+\kappa|z|+\psi(y)|z|^2$, where $\psi(y)\in C_+(D)$ (can be singular at 0 when $D=(0,\infty)$), $\delta\in H_1^1$, $\gamma\in\mathcal{C}$ are nonnegative, and $c>0, \kappa\geq0$ are constants. Then $g$ satisfies the assumption in Proposition \ref{pr5.2} with $f(y):=1_{\{y<c\}}\psi(c)+1_{\{y\geq c\}}\sup_{x\in[c,y]}\psi(x).$

(ii) In Proposition \ref{pr5.2}, since $f$ is nondecreasing, by Example \ref{exC.2}(i), when $f(y)\geq
\beta>0$ for some $y\in D$, to guarantee the existence of a solution, $\xi$ satisfies at least the integrability: $\exp(2\beta\Lambda^{\delta,\gamma}_T(\xi))\in L^p({\mathcal{F}}_T)$. In fact, from the proof of Proposition \ref{pr5.2}, it can be seen that the assumption that $f$ is nondecreasing is used to obtain the crucial inequalities: ``$g_2(\cdot,\cdot\vee\Lambda_t^{\delta,\gamma}(L_t\vee c),\cdot)\leq g_1(\cdot,\cdot\vee\Lambda_t^{\delta,\gamma}(L_t\vee c),\cdot)$" in Step 2, and ``$g_4(\cdot,\cdot,\cdot)\leq g_3(\cdot,\cdot,\cdot)$" in Step 3. One can check that the two inequalities also hold in the following case:
\begin{itemize}
       \item $\delta_t\equiv0$ and $D=(0,\infty)$, $f(y)=\frac{\beta}{y^r}$, $0< r\leq1$, $\beta>0$.
\end{itemize}
Hence, the integrability of $\xi$ in Proposition \ref{pr5.2} may be weakened in some cases.
\end{remark}

In view of Remark \ref{r5.4}(ii), we consider the following assumption.
\begin{itemize}
  \item \textbf{(5A1')}\label{5A1'} $D=(0,\infty)$ and there exist nonnegative processes $\gamma, \vartheta\in{\mathcal{C}}$ and constants $\kappa\geq0,$ $\beta\geq0,$ $\nu\geq0,$ $c\geq b>0,$ such that $dt\times dP\textrm{-}a.e.,$ for each $(y,z)\in D\times{\mathbf{ R}}^{\mathit{d}},$
      \begin{itemize}
  \item\textbf{(i)} \ \ \ \ \ \ \ \ \ \ \ \ \ \ \ \ \ \ \ \ \ \ \ \ \ \ \ \  $1_{\{y\geq c\}}{g}(t,y,z)\leq\gamma_t |y|+\kappa|z|+\frac{\beta}{y}|z|^2;$
  \item\textbf{(ii)} \ \ \ \ \ \ \ \ \ \ \ \ \ \ \ \ \ \ \ \ \ \ \ \ \ \ $1_{\{y\leq b\}}{g}(t,y,z)\geq-\vartheta_t|y|-\kappa|z|-\frac{\nu}{y}|z|^2.$
\end{itemize}
\end{itemize}

\begin{cor}\label{c5.5} Let (\hyperref[5A1']{5A1'}-(i)) and (\hyperref[4A2']{4A2'}) hold. If $|e^{\int_0^T\gamma_sds}(\xi\vee c)|^{1+2\beta}\in L^p({\mathcal{F}}_T)$ and $|e^{\int_0^\cdot\gamma_sds}L|^{1+2\beta}\in{\mathcal{S}}^p$, then the RBSDE$(g,\xi,L)$ admits a unique solution $(Y,Z,K)$ such that $|e^{\int_0^\cdot\gamma_sds}Y|^{1+2\beta}\in{\mathcal{S}}^p$.
\end{cor}
\begin{proof} \emph{Existence:} By Example \ref{exC.2}(ii), we have $u_{\frac{\beta}{\cdot}}(x)=\frac{1}{1+2\beta}(x^{1+2\beta}-1)$. \footnote{The notation $u_{\frac{\beta}{\cdot}}(x)$ can be found in Example \ref{exC.2} in Appendix \ref{appC}.} This implies $u_{\frac{\beta}{\cdot}}(\Lambda^{0, \gamma}(\xi\vee c))\in L^p({\mathcal{F}}_T)$ and $u_{\frac{\beta}{\cdot}}(\Lambda^{0, \gamma}(L))\in{\mathcal{S}}^p$. Then, by Proposition \ref{pr5.2}, Remark \ref{r5.4}(ii), and (\hyperref[5A1']{5A1'}-(i)), we deduce that the RBSDE$(g,\xi,L)$ admits a minimal solution $(Y,Z,K)$ such that $|e^{\int_0^\cdot\gamma_sds}Y|^{1+2\beta}\in{\mathcal{S}}^p$, and a maximal solution $(\bar{Y},\bar{Z},\bar{K})$ such that $|e^{\int_0^\cdot\gamma_sds}\bar{Y}|^{1+2\beta}\in{\mathcal{S}}^p$.

\emph{Uniqueness:} Let $(\hat{Y},\hat{Z},\hat{K})$ be another solution to the RBSDE$(g,\xi,L)$ such that $|e^{\int_0^\cdot\gamma_sds}\hat{Y}|^{1+2\beta}\in{\mathcal{S}}^p$. By (\hyperref[4A2']{4A2'}) and Remark \ref{r4.1}(i), we have $dt\times dP\textrm{-}a.e.,$
\begin{equation*}
\frac{1}{{1-\theta}}(g(t,\hat{Y}_t,\hat{Z}_t)-\theta g(t,{Y}_t,{Z}_t))\leq g\left(t,\frac{\hat{Y}_t-\theta {Y}_t}{1-\theta},\frac{\hat{Z}_t-\theta {Z}_t}{1-\theta}\right).
\end{equation*}
Moreover, for each $\theta\in(0,1)$, we have $u_{\frac{\beta}{\cdot}}(\Lambda^{0,\gamma}(\frac{\hat{Y}-\theta Y}{1-\theta}))\in{\mathcal{S}}^p$. Then, by the Existence above, we get that for each $\theta\in(0,1)$, the RBSDE$(g,\xi,\frac{\hat{Y}-\theta Y}{1-\theta})$ admits a solution $(y^\theta,z^\theta,k^\theta)$ such that $|e^{\int_0^\cdot\gamma_sds}y^\theta|^{1+2\beta}\in{\mathcal{S}}^p$, and the RBSDE$(g,\xi,\hat{Y}_t)$ admits a maximal solution $(\overline{Y},\overline{Z},\overline{K})$ such that $|e^{\int_0^\cdot\gamma_sds}\overline{Y}|^{1+2\beta}\in{\mathcal{S}}^p$.
Since $Y_t\leq \hat{Y}_t$, by a argument as in Step 3 of the proof of Theorem \ref{th4.2} (replacing $G$ therein with $g$), we deduce that for each $\theta\in(0,1)$, $(y^\theta_t,z^\theta_t,k^\theta_t)$ is a solution to the RBSDE$(g,\xi,\hat{Y}_t)$. Thus, for each $\theta\in(0,1)$, $\overline{Y}_t\geq\frac{\hat{Y}_t-\theta Y_t}{1-\theta}.$ This implies $(\hat{Y}_t,\hat{Z}_t,\hat{K}_t)=(Y_t,Z_t,K_t).$
\end{proof}
\subsection{Quadratic BSDEs}
When $D=(0,\infty)$, we have the following well-posedness results for unbounded solutions of BSDEs.
\begin{proposition}\label{pr5.6} Let $D=(0,\infty)$, and let (\hyperref[5A1]{5A1}-(i)(iii)) hold with $f(y)$ nondecreasing. Let $u_{\frac{-\nu}{\cdot}}(e^{-\int_0^T\vartheta_sds}(\xi\wedge b)),$ $u_f(\Lambda^{\delta,\gamma}_T(\xi\vee c))\in L^p({\mathcal{F}}_T)$. Then the BSDE$(g,\xi)$ admits at least one solution $(Y,Z)$ such that $u_{\frac{-\nu}{\cdot}}(e^{-\int_0^\cdot\vartheta_sds}Y)$, $u_f(\Lambda^{\delta,\gamma}(Y))\in{\mathcal{S}}^p$. Moreover, if we further assume that for each $q\geq1$, $u_{f}(q\Lambda^{\delta, \gamma}_T(\xi\vee c))\in L^p({\mathcal{F}}_T)$ and one of the following two conditions holds:

(i) $\int_0^T\gamma_tdt\in L^\infty(\mathcal{F}_T)$ and (\hyperref[5A2]{5A2}) holds with $\int_0^T(\tilde{\mu}_t+\tilde{\gamma}_t)dt\in L^\infty(\mathcal{F}_T)$;

(ii) (\hyperref[4A2']{4A2'}) holds,
\\then the BSDE$(g,\xi)$ admits a unique solution $(Y,Z)$ such that $u_{\frac{-\nu}{\cdot}}(e^{-\int_0^\cdot\vartheta_sds}Y)\in{\mathcal{S}}^p$ and for each $q\geq1$, $u_{f}(q\Lambda^{\delta, \gamma}(Y))\in{\mathcal{S}}^p$, in particular,
\begin{equation*}\label{5.9}
\left\{
  \begin{array}{ll}
   |e^{-\int_0^\cdot\vartheta_sds}Y|^{1-2\nu}\in{\mathcal{S}}^p\ and\ \exp(\Lambda^{\delta,\gamma}(Y))\in\bigcap_{q\geq1}{\mathcal{S}}^q, & \textrm{when}\ \nu\neq\frac{1}{2},\ f\equiv\beta>0;\\
   |e^{-\int_0^\cdot\vartheta_sds}Y|^{1-2\nu},\ \Lambda^{\delta,\gamma}(Y)\in{\mathcal{S}}^p, & \textrm{when}\ \nu\neq\frac{1}{2},\ f\equiv0;\\
\ln(e^{-\int_0^\cdot\vartheta_sds}Y)\in{\mathcal{S}}^p\ and\ \exp(\Lambda^{\delta,\gamma}(Y))\in\bigcap_{q\geq1}{\mathcal{S}}^q, & \textrm{when}\ \nu=\frac{1}{2},\ f\equiv\beta>0;\\
\ln(e^{-\int_0^\cdot\vartheta_sds}Y),\ \Lambda^{\delta,\gamma}(Y)\in{\mathcal{S}}^p, & \textrm{when}\ \nu=\frac{1}{2},\ f\equiv0.
  \end{array}
\right.\tag{5.9}
\end{equation*}
\end{proposition}

\begin{proof} \emph{Existence:} Since $u_{\frac{-\nu}{\cdot}}(e^{-\int_0^T\vartheta_sds}(\xi\wedge b))\in L^p({\mathcal{F}}_T)$ and for all $t\in[0,T]$,
\begin{equation*}
  u_{\frac{-\nu}{\cdot}}(e^{-\int_0^T\vartheta_sds}(\xi\wedge b))\leq u_{\frac{-\nu}{\cdot}}(e^{-\int_0^t\vartheta_sds}b)\leq u_{\frac{-\nu}{\cdot}}(b),
\end{equation*}
we have $u_{\frac{-\nu}{\cdot}}(e^{-\int_0^\cdot\vartheta_sds}b)\in{\mathcal{S}}^p$. Then, by Lemma \ref{lmE.1}, the RBSDE$(-\vartheta_ty+\kappa|z|-\frac{\nu}{y}|z|^2,-(\xi\wedge b),-b)$ admits a unique solution $(y^1,z^1,k^1)$ such that $-y^1\in\mathcal{C}_D$ and $u_{\frac{-\nu}{\cdot}}(-e^{-\int_0^\cdot\vartheta_sds}y^1)\in{\mathcal{S}}^p$. Set $(Y^1_t,Z^1_t):=(-y^1_t,-z^1_t)$, we get that $0<Y_t^1\leq b$ and
\begin{equation*}\label{5.10}
Y^1_t=\xi\wedge b+\int_t^T(-\vartheta_s|Y^1_s|-\kappa|Z^1_s|-\frac{\nu}{Y^1_s}|Z^1_s|^2)ds-k^1_T+k^1_t-\int_t^TZ^1_s\cdot dB_s,\quad t\in[0,T].\tag{5.10}
\end{equation*}

Since $u_{\frac{-\nu}{\cdot}}(e^{-\int_0^T\vartheta_sds}(\xi\wedge b)),$ $u_f(\Lambda^{\delta,\gamma}_T(\xi\vee c))\in L^p({\mathcal{F}}_T)$ and for all $t\in[0,T]$,
\begin{equation*}
  u_{\frac{-\nu}{\cdot}}(e^{-\int_0^T\vartheta_sds}(\xi\wedge b))\leq u_f(\Lambda^{\delta,\gamma}_t(c))\leq u_f(\Lambda^{\delta,\gamma}_T(\xi\vee c)),
\end{equation*}
we have $u_f(\Lambda^{\delta,\gamma}_t(c))\in{\mathcal{S}}^p$. Then by Proposition \ref{pr5.2}, the RBSDE$(\delta_t+\gamma_t|y|+\kappa|z|+f(y)|z|^2,\xi\vee c, c)$ admits a solution $(Y^2,Z^2,K^2)$ such that $u_f(\Lambda^{\delta,\gamma}(Y^2))\in{\mathcal{S}}^p$. In view of $0<Y_t^1\leq b\leq c\leq Y^2_t$, by (\hyperref[5A1]{5A1}-(i)(iii)), (\ref{5.10}) and Lemma \ref{l2.4}, the BSDE$(g,\xi)$ admits a minimal solution $(Y,Z)$ such that $Y_t\geq Y_t^1$. Moreover, we also have $Y_t^1\leq Y_t\leq Y_t^2$, which implies
\begin{eqnarray*}
u_{\frac{-\nu}{\cdot}}(e^{-\int_0^t\vartheta_sds}Y^1_t)\leq u_{\frac{-\nu}{\cdot}}(e^{-\int_0^t\vartheta_sds}Y_t)\leq u_f(e^{-\int_0^t\vartheta_sds}Y_t)\leq u_f(\Lambda^{\delta,\gamma}_t(Y_t))\leq u_f(\Lambda^{\delta,\gamma}_t(Y^2_t)).
\end{eqnarray*}
It follows that $u_{\frac{-\nu}{\cdot}}(e^{-\int_0^\cdot\vartheta_sds}Y)$, $u_f(\Lambda^{\delta,\gamma}(Y))\in{\mathcal{S}}^p$.

\emph{Uniqueness:} Since $u_{\frac{-\nu}{\cdot}}(e^{-\int_0^T\vartheta_sds}(\xi\wedge b))\in L^p({\mathcal{F}}_T)$ and for each $q\geq1$, $u_{f}(q\Lambda^{\delta, \gamma}_T(\xi\vee c))\in L^p({\mathcal{F}}_T)$, by (\ref{5.5}) and the proof of existence above, we deduce that $(Y,Z)$ is a minimal solution to the BSDE$(g,\xi)$ such that $Y_t\geq Y_t^1$, and $u_{\frac{-\nu}{\cdot}}(e^{-\int_0^\cdot\vartheta_sds}Y)\in{\mathcal{S}}^p$ and for each $q\geq1$, $u_{f}(q\Lambda^{\delta,\gamma}(Y))\in{\mathcal{S}}^p$.

Let $(\hat{Y},\hat{Z})$ be another solution to the BSDE$(g,\xi)$ such that $u_{\frac{-\nu}{\cdot}}(e^{-\int_0^\cdot\vartheta_sds}\hat{Y})\in{\mathcal{S}}^p$ and for each $q\geq1$, $u_{f}(q\Lambda^{\delta,\gamma}(\hat{Y}))\in{\mathcal{S}}^p$.
In the following, we will prove that the BSDE$(g,\xi)$ admits a solution $(\tilde{y},\tilde{z})$ satisfying $u_{\frac{-\nu}{\cdot}}(e^{-\int_0^\cdot\vartheta_sds}\tilde{y})\in{\mathcal{S}}^p$ and for each $q\geq1$, $u_{f}(q\Lambda^{\delta,\gamma}(\tilde{y}))\in{\mathcal{S}}^p$ such that $\tilde{y}_t\leq \hat{Y}_t$ and $\tilde{y}_t\leq Y_t$. This is a crucial step of this proof.

Since
\begin{equation*}
u_{\frac{-\nu}{\cdot}}(e^{-\int_0^t\vartheta_sds}(\hat{Y}_t\wedge Y_t^1))= u_{\frac{-\nu}{\cdot}}(e^{-\int_0^t\vartheta_sds}(\hat{Y}_t))\wedge u_{\frac{-\nu}{\cdot}}(e^{-\int_0^t\vartheta_sds}(Y_t^1)),
\end{equation*}
we have $u_{\frac{-\nu}{\cdot}}(e^{-\int_0^\cdot\vartheta_sds}(\hat{Y}\wedge Y^1))\in{\mathcal{S}}^p.$ By Lemma \ref{lmE.1} again, the RBSDE$(-\vartheta_ty+\kappa|z|-\frac{\nu}{y}|z|^2,-(\xi\wedge b),-(\hat{Y}\wedge Y^1))$ admits a unique solution $(\tilde{y}^1,\tilde{z}^1,\tilde{k}^1)$ such that $-\tilde{y}^1\in\mathcal{C}_D$ and $u_{\frac{-\nu}{\cdot}}(-e^{-\int_0^\cdot\vartheta_sds}\tilde{y}^1)\in{\mathcal{S}}^p$. Set $(\tilde{Y}^1_t,\tilde{Z}^1_t):=(-\tilde{y}^1_t,-\tilde{z}^1_t)$, we have
\begin{equation*}\label{5.11}
\tilde{Y}^1_t=\xi\wedge b+\int_t^T(-\vartheta_s|\tilde{Y}^1_s|-\kappa|\tilde{Z}^1_s|-\frac{\nu}{\tilde{Y}^1_s}|\tilde{Z}^1_s|^2)ds-
\tilde{k}^1_T+\tilde{k}^1_t-\int_t^T\tilde{Z}^1_s\cdot dB_s,\quad t\in[0,T].\tag{5.11}
\end{equation*}
Since $0<\tilde{Y}^1_t\leq \hat{Y}_t\wedge Y_t^1\leq b$, by (\hyperref[5A1]{5A1}-(iii)) and (\ref{5.11}), it follows that the BSDE$(g,\xi)$ is dominated by $\tilde{Y}^1$ and $\hat{Y}_t$. Then by Lemma \ref{l2.4}, the BSDE$(g,\xi)$ admits a minimal solution $(\tilde{y},\tilde{z})$ such that $\tilde{y}_t\geq\tilde{Y}^1_t$. Moreover, since $0<\tilde{Y}^1_t\leq \tilde{y}_t\leq \hat{Y}_t$, we have for each $q\geq1$,
\begin{eqnarray*}
u_{\frac{-\nu}{\cdot}}(e^{-\int_0^t\vartheta_sds}\tilde{Y}^1_t)\leq u_{\frac{-\nu}{\cdot}}(e^{-\int_0^t\vartheta_sds}\tilde{y}_t)\leq u_{f}(e^{-\int_0^t\vartheta_sds}\tilde{y}_t)\leq u_{f}(q\Lambda^{\delta,\gamma}_t(\tilde{y}_t))
\leq u_{f}(q\Lambda^{\delta,\gamma}_t(\hat{Y}_t)).
\end{eqnarray*}
It follows that $u_{\frac{-\nu}{\cdot}}(e^{-\int_0^\cdot\vartheta_sds}\tilde{y})\in{\mathcal{S}}^p$ and for each $q\geq1$, $u_{f}(q\Lambda^{\delta,\gamma}(\tilde{y}))\in{\mathcal{S}}^p$.

Since $(\tilde{y},\tilde{z})$ (resp. $(Y,Z)$) is a minimal solution to the BSDE$(g,\xi)$ such that $\tilde{y}_t\geq\tilde{Y}^1_t$ (resp. $Y_t\geq Y_t^1$) and $\tilde{Y}^1_t\leq Y_t^1$, we have $\tilde{y}_t\leq Y_t$. Thus, $(\tilde{y},\tilde{z})$ is a solution to the BSDE$(g,\xi)$ satisfying $u_{\frac{-\nu}{\cdot}}(e^{-\int_0^\cdot\vartheta_sds}\tilde{y})\in{\mathcal{S}}^p$ and for each $q\geq1$, $u_{f}(q\Lambda^{\delta,\gamma}(\tilde{y}))\in{\mathcal{S}}^p$, such that $\tilde{y}_t\leq \hat{Y}_t$ and $\tilde{y}_t\leq Y_t$.

\textbf{Case (i):} By (\hyperref[5A2]{5A2}) and (\hyperref[5A1]{5A1}-(i)), we have $dt\times dP\textrm{-}a.e.,$ for each $(y,z)\in D\times{\mathbf{R}}^{\mathit{d}},$
\begin{equation*}
1_{\{y\geq c\}}G(t,y,z,\tilde{y}_t,Y_t)\leq\bar{\delta}_t+\bar{\gamma}_t|y|
+\bar{\kappa}|z|+f(y)|z|^2,
\end{equation*}
where $\bar{\delta}_t:=\tilde{\mu}_t(1+|\tilde{y}_t|+|Y_t|)+\delta_t$, $\bar{\gamma}_t:=\tilde{\gamma}_t+\gamma_t$ and $\bar{\kappa}=\tilde{\kappa}+\kappa.$ As discussed in (\ref{5.6})-(\ref{5.8}), we can deduce that $u_{f}(\Lambda^{\bar{\delta},\bar{\gamma}}_T(\xi\vee c))\in L^p({\mathcal{F}}_T)$ and
$u_{f}(\Lambda^{\bar{\delta},\bar{\gamma}}(Y)), u_{f}(\Lambda^{\bar{\delta},\bar{\gamma}}(\frac{Y-\theta\tilde{y}}{1-\theta}))\in{\mathcal{S}}^p$ for each $\theta\in(0,1)$.
Then, since $\tilde{y}_t\leq Y_t$, by a similar argument as in Case (i) of the proof of Proposition \ref{pr5.3} (consider the $\theta$-difference process $\frac{Y-\theta\tilde{y}}{1-\theta}$, correspondingly), we obtain $(\tilde{y}_t,\tilde{z}_t)=({Y}_t,{Z}_t).$ Similarly, since $\tilde{y}_t\leq \hat{Y}_t,$ we can also get $(\tilde{y}_t,\tilde{z}_t)=(\hat{Y}_t,\hat{Z}_t).$ Thus, we have $(Y_t,Z_t)=(\hat{Y}_t,\hat{Z}_t),$ which implies that $(Y_t,Z_t)$ is a unique solution to the BSDE$(g,\xi)$ such that $u_{\frac{-\nu}{\cdot}}(e^{-\int_0^\cdot\vartheta_sds}Y)\in{\mathcal{S}}^p$ and for each $q\geq1$, $u_{f}(q\Lambda^{\delta, \gamma}(Y))\in{\mathcal{S}}^p$.

\textbf{Case (ii):} The proof is similar to Case (i) (replacing $G$ therein with $g$). We omit it.

By Example \ref{exC.2}(i)(ii), we further obtain (\ref{5.9}). The proof is complete.
\end{proof}

\begin{cor}\label{c5.7} Let (\hyperref[5A1']{5A1'}-(i)(ii)) and (\hyperref[4A2']{4A2'}) hold. Then the following hold:

(i) If $\nu\neq\frac{1}{2}$ and $|e^{-\int_0^T\vartheta_sds}(\xi\wedge b)|^{1-2\nu}$, $|e^{\int_0^T\gamma_sds}(\xi\vee c)|^{1+2\beta}\in L^p({\mathcal{F}}_T)$, then the BSDE$(g,\xi)$ admits a unique solution $(Y,Z)$ such that $|e^{-\int_0^\cdot\vartheta_sds}Y|^{1-2\nu}$, $|e^{\int_0^\cdot\gamma_sds}Y|^{1+2\beta}\in{\mathcal{S}}^p$.

(ii) If $\nu=\frac{1}{2}$ and $\ln(e^{-\int_0^T\vartheta_sds}(\xi\wedge b))$, $|e^{\int_0^T\gamma_sds}(\xi\vee c)|^{1+2\beta}\in L^p({\mathcal{F}}_T)$, then the BSDE$(g,\xi)$ admits a unique solution $(Y,Z)$ such that $\ln(e^{-\int_0^\cdot\vartheta_sds}Y)$, $|e^{\int_0^\cdot\gamma_sds}Y|^{1+2\beta}\in{\mathcal{S}}^p$.
\end{cor}
\begin{proof} We only consider (i), (ii) is similar. By Lemma \ref{lmE.1}, the RBSDE$(-\vartheta_ty+\kappa|z|-\frac{\nu}{y}|z|^2,-(\xi\wedge b),-b)$ admits a unique solution $(y^1,z^1,k^1)$ such that $-y^1\in\mathcal{C}_D$ and $|e^{-\int_0^\cdot\vartheta_sds}y^1|^{1-2\nu}\in{\mathcal{S}}^p$. By Corollary \ref{c5.5}, the RBSDE$(\gamma_t|y|+\kappa|z|+\frac{\beta}{y}|z|^2,\xi\vee c, c)$ admits a unique solution $(Y^2,Z^2,K^2)$ such that $|e^{\int_0^\cdot\gamma_sds}Y^2|^{1+2\beta}\in{\mathcal{S}}^p$. Then, by a similar argument as in the existence in the proof of Proposition \ref{pr5.6}, we deduce that the BSDE$(g,\xi)$ admits a minimal solution $(Y,Z)$ such that $Y_t\geq-y^1$, and $|e^{-\int_0^\cdot\vartheta_sds}Y|^{1-2\nu}$, $|e^{\int_0^\cdot\gamma_sds}Y|^{1+2\beta}\in{\mathcal{S}}^p$.

Let $(\hat{Y},\hat{Z})$ be another solution to the BSDE$(g,\xi)$ such that $|e^{-\int_0^\cdot\vartheta_sds}\hat{Y}|^{1-2\nu}$, $|e^{\int_0^\cdot\gamma_sds}\hat{Y}|^{1+2\beta}\in{\mathcal{S}}^p$. By a similar argument as in the uniqueness in the proof of Proposition \ref{pr5.6}, we can get that the BSDE$(g,\xi)$ admits a solution $(\tilde{y},\tilde{z})$ satisfying $|e^{-\int_0^\cdot\vartheta_sds}Y|^{1-2\nu}$, $|e^{\int_0^\cdot\gamma_sds}Y|^{1+2\beta}\in{\mathcal{S}}^p$, such that $\tilde{y}_t\leq \hat{Y}_t$ and $\tilde{y}_t\leq Y_t$. Then, by a similar argument as in the uniqueness in the proof of Corollary \ref{c5.5}, we can get that $\tilde{y}_t=\hat{Y}_t$ and $\tilde{y}_t=Y_t$. Thus, $(Y,Z)$ is a unique solution to the BSDE$(g,\xi)$ such that $|e^{-\int_0^\cdot\vartheta_sds}Y|^{1-2\nu}$, $|e^{\int_0^\cdot\gamma_sds}Y|^{1+2\beta}\in{\mathcal{S}}^p$.
\end{proof}

\begin{remark}\label{r5.8} Recently, an existence and uniqueness result was obtained by \cite[Theorem 17]{LG} for a generator $g$ satisfying (\hyperref[4A2']{4A2'}) and the following condition:
\begin{equation*}\label{5.12}
0\leq g(t,y,z)\leq\gamma_ty+ \alpha_ty|\ln(y)|+\kappa_t\cdot z+\frac{\beta}{y}|z|^2,\ \  (y,z)\in(0,\infty)\times\mathbf{R}^d,\tag{5.12}
\end{equation*}
where $\gamma, \alpha\in H_1^2,\kappa\in H_d^2$ are all positive and bounded, and $\beta>0$ is a constant. In comparison, the process $\gamma$ in Corollary \ref{c5.7} is unbounded, and the generator $g$ in Corollary \ref{c5.7} can grow arbitrarily fast to $+\infty$ as $y\to0^+$ and can grow arbitrarily fast to $-\infty$ or grow linearly to $+\infty$ as $y\to+\infty$, whereas the generator $g$ in (\ref{5.12}) exhibits a superlinear growth to $+\infty$ as $y\to+\infty$. Moreover, \cite[Theorem 17]{LG} requires a stronger integrability condition on $\xi$ than that in Corollary \ref{c5.7} even when $\alpha\equiv0$ in (\ref{5.12}). We also refer to \cite[Theorem 2.3]{BT} for a uniqueness result on bounded solutions to a BSDE with generator $g$ satisfying (\hyperref[4A2']{4A2'}) and $0\leq g(t,y,z)\leq\delta_t+\gamma_ty+\kappa_t\cdot z+\frac{\beta}{y}|z|^2, (y,z)\in(0,\infty)\times\mathbf{R}^d$, where $\delta, \gamma\in\mathcal{S}^\infty$ and $\kappa\in \mathcal{H}_d^{BMO}$ are all positive, and $\beta>0$ is a constant.
\end{remark}

When $D=\mathbf{R}$, we have the following well-posedness results for unbounded solutions of BSDEs, which can be seen as an extension of the corresponding results in \cite{BH08, B17, B19, FHT2}.
\begin{proposition}\label{pr5.9} Let $D=\mathbf{R}$, and let  (\hyperref[5A1]{5A1}-(i)(ii)) hold with $f(y)$ nondecreasing. Let $u_f(\Lambda^{\delta,\gamma}_T(c)),$ $u_f(\Lambda^{\delta,\gamma}_T(\xi))\in L^p({\mathcal{F}}_T)$. Then the BSDE$(g,\xi)$ admits at least one solution $(Y,Z)$ such that $u_f(\Lambda^{\delta,\gamma}(Y))\in{\mathcal{S}}^p.$ Moreover, if we further assume that for each $q\geq1$, $u_f(q\Lambda^{\delta,\gamma}_T(c)),$ $u_f(q\Lambda^{\delta,\gamma}_T(\xi))\in L^p({\mathcal{F}}_T)$ and one of the following two conditions holds:

(i) $\int_0^T\gamma_tdt\in L^\infty(\mathcal{F}_T)$ and (\hyperref[5A2]{5A2}) holds with $\int_0^T(\tilde{\mu}_t+\tilde{\gamma}_t)dt\in L^\infty(\mathcal{F}_T)$;

(ii) (\hyperref[4A2']{4A2'}) holds,\\
then the BSDE$(g,\xi)$ admits a unique solution $(Y,Z)$ such that for each $q\geq1$, $u_{f}(q\Lambda^{{\delta},{\gamma}}(Y))\in{\mathcal{S}}^p$, in particular
\begin{equation*}\label{5.13}
\left\{
  \begin{array}{ll}
    \exp(\Lambda^{\delta,\gamma}(Y))\in\bigcap_{q\geq1}{\mathcal{S}}^q,\ & when\ f\equiv\beta>0;\\
    \Lambda^{\delta,\gamma}(Y)\in{\mathcal{S}}^p, & when\ f\equiv0.
  \end{array}
\right.\tag{5.13}
\end{equation*}
\end{proposition}
\begin{proof}\emph{Existence:} Since $u_f(\Lambda^{\delta,\gamma}_T(c)), u_f(\Lambda^{\delta,\gamma}_T(\xi))\in L^p({\mathcal{F}}_T)$, we have $u_f(\Lambda^{\delta,\gamma}_T(\xi\vee c)), u_f(\Lambda^{\delta,\gamma}_T((-\xi)\vee c))\in L^p({\mathcal{F}}_T)$. By Proposition \ref{pr5.2}, the RBSDE$(\delta_t+\gamma_t |y|+\kappa|z|+f(y)|z|^2,\xi\vee c,c)$ admits a minimal solution $(Y^1,Z^1,K^1)$ such that $Y^1\in\mathcal{C}$, in particular, $u_f(\Lambda^{\delta,\gamma}_T(Y^1))\in{\mathcal{S}}^p$, and the RBSDE$(\delta_t+\gamma_t|y|+\kappa|z|+f(y)|z|^2,(-\xi)\vee c,c)$ admits a minimal solution $(Y^2,Z^2,K^2)$ such that $Y^2\in\mathcal{C}$, in particular, $u_f(\Lambda^{\delta,\gamma}_T(Y^2))\in{\mathcal{S}}^p$.
By (\hyperref[5A1]{5A1}-(i)(ii)) and the fact that $-Y^2_t\leq -c\leq c\leq Y^1_t$, we deduce that the BSDE$(g,\xi)$ is dominated by $-Y^2$ and $Y^1$. Lemma \ref{l2.4} then implies that the BSDE$(g,\xi)$ admits a minimal solution $(Y,Z)$ such that $Y_t\geq -Y^2_t.$ Moreover, we have $u_f(\Lambda^{\delta,\gamma}_T(Y))\in{\mathcal{S}}^p.$

\emph{Uniqueness:} Since for each $q\geq1$, $u_f(q\Lambda^{\delta,\gamma}_T(c)),$ $u_f(q\Lambda^{\delta,\gamma}_T(\xi))\in L^p({\mathcal{F}}_T)$, by Proposition \ref{pr5.2}(i), we can deduce that for each $q\geq1$, $u_{f}(q\Lambda^{{\delta},{\gamma}}(Y^i))\in{\mathcal{S}}^p$, $i=1,2$. This implies that $(Y,Z)$ is a solution to the BSDE$(g,\xi)$ such that for each $q\geq1$, $u_{f}(q\Lambda^{{\delta},{\gamma}}(Y))\in{\mathcal{S}}^p.$

We only consider Case (i), Case (ii) is similar. Let the condition (i) hold. Let $(\hat{Y},\hat{Z})$ be another solution to the BSDE$(g,\xi)$ such that for each $q\geq1$, $u_{f}(\Lambda^{\delta,\gamma}(\hat{Y}))\in{\mathcal{S}}^p$. Since (\hyperref[5A2]{5A2}) and (\hyperref[5A1]{5A1}-(i)) holds, we have $dt\times dP\textrm{-}a.e.,$ for each $(y,z)\in D\times{\mathbf{R}}^{\mathit{d}},$
\begin{equation*}
1_{\{y\geq \tilde{c}\}}G(t,y,z,\hat{Y}_t,Y_t)\leq\bar{\delta}_t+\bar{\gamma}_t|y|
+\bar{\kappa}|z|+f(y)|z|^2,
\end{equation*}
where $\bar{\delta}_t:=\tilde{\mu}_t(1+|\hat{Y}|+|Y_t|)+\delta_t$, $\bar{\gamma}_t:=\tilde{\gamma}_t+\gamma_t$ and $\bar{\kappa}=\tilde{\kappa}+\kappa.$ As discussed in (\ref{5.6})-(\ref{5.8}), we can deduce that $u_{f}(\Lambda^{\bar{\delta},\bar{\gamma}}_T(\xi\vee c))\in L^p({\mathcal{F}}_T)$ and
$u_{f}(\Lambda^{\bar{\delta},\bar{\gamma}}(\hat{Y})), u_{f}(\Lambda^{\bar{\delta},\bar{\gamma}}(\frac{\hat{Y}-\theta Y}{1-\theta}))\in{\mathcal{S}}^p$ for each $\theta\in(0,1)$. By Proposition \ref{pr5.2}, for each $\theta\in(0,1)$, the RBSDE$(G(t,y,z,\hat{Y}_t,Y_t),\xi,\frac{\hat{Y}-\theta Y}{1-\theta})$ admits a minimal solution $(y^\theta,z^\theta,k^\theta)$ such that $u_{f}(\Lambda^{\bar{\delta},\bar{\gamma}}(y^\theta))\in{\mathcal{S}}^p$. Then, by a similar argument as in Case (i) of the proof of Theorem \ref{th4.3}, and Proposition \ref{pr5.2}, we can deduce that the BSDE$(G(t,y,z,\hat{Y}_t,Y_t),\xi)$ admits a maximal solution $(\hat{y}_t,\hat{z}_t)$ such that for each $\theta\in(0,1)$, $\hat{y}_t\geq\frac{\hat{Y}_t-\theta Y_t}{1-\theta}$, which gives $\hat{Y}_t\leq Y_t$. Similarly, by considering the difference $\frac{Y-\theta\hat{Y}}{1-\theta}$, we can also get $\hat{Y}_t\geq Y_t$. Thus $(\hat{Y}_t,\hat{Z}_t)=(Y_t,Z_t)$, which implies that $(Y_t,Z_t)$ is a unique solution to the BSDE$(g,\xi)$ such that for each $q\geq1$, $u_f(q\Lambda^{\delta,\gamma}(Y_t))\in{\mathcal{S}}^p$.

By Example \ref{exC.2}(i), we further get (\ref{5.13}). The proof is complete.
\end{proof}

A key difference between the quadratic growth conditions in (\hyperref[5A1]{5A1}) (resp. (\hyperref[4A1]{4A1})) and those in existing
studies lies in the one-sided growth in $y$. These one-sided growth conditions contain singular
generators and generators with general stochastic coefficients. Some examples of (\hyperref[4A1]{4A1}) were given in Example \ref{ex4.6}, we now show some examples of (\hyperref[5A1]{5A1}).

\begin{example}\label{ex5.10}
(i) Let $D={\mathbf{R}}$ and
\begin{equation*}
\mathfrak{g}(t,y,z)=\delta_t+k(y)(\phi(t,y)+\psi(t,y)|z|^r)+\gamma_t|y|+f(y)|z|^2,
\end{equation*}
where $\delta\in H^1_1$, $\gamma\in{\mathcal{C}}$, $k(y)$, $f(y)\in C({\mathbf{R}})$, $0<r\leq2$, and $\phi(\omega,t,y),\psi(\omega,t,y): \Omega\times [0,T]\times{\mathbf{R}}\longmapsto{\mathbf{R}}$ are both measurable with respect to ${\mathcal{P}}\otimes{\mathcal{B}}({\mathbf{R}})$ and continuous on $[0,T]\times {\mathbf{R}}$.
By Remark \ref{r2.1}, $\mathfrak{g}$ satisfies (\hyperref[2A1]{2A1}). If there exists a constant $c>0$ such that for each $y\geq c$, $k(y)=0$, then $\mathfrak{g}$ satisfies (\hyperref[5A1]{5A1}-(i)); and moreover, if we also have for each $y\leq -c$, $k(y)=0$, then $\mathfrak{g}$ satisfies (\hyperref[5A1]{5A1}-(i)(ii)).

(ii) Let $D=(0,\infty)$ and
\begin{equation*}
\mathfrak{g}(t,y,z)=\delta_t+\frac{b_1}{y^k}+b_2|y|+b_3\phi(t,y)+b_4|z|^l+\frac{b_5}{y^r}|z|^2+g_1(z),
\end{equation*}
where $\delta\in H^1_1$ is nonnegative, $k>0$, $1\leq l\leq2$ and $b_i\geq0$ $(1\leq i\leq5)$ are constants,  $\phi(t,y)$ is the function defined in (\ref{4.10}), and $g_1(z):{\bf{R}}^d\rightarrow{\bf{R}}$ is a bounded Lipschitz function with a bounded support and $g(0)\geq0$. By Remark \ref{r2.1}, Remark \ref{r4.1}(i)(iii) and Lemma \ref{lmD.1}, we can check that $\mathfrak{g}$ satisfies (\hyperref[2A1]{2A1}), (\hyperref[5A1]{5A1}-(i)(iii)), and (\hyperref[5A2]{5A2}).
\end{example}

\begin{remark}
All the uniqueness results in Theorems \ref{th4.2} and \ref{th4.3}, and Propositions \ref{pr5.3}, \ref{pr5.6} and \ref{pr5.9} hold true for any continuous nondecreasing function $f$, however,
up to now, we do not find that there exists a function $g$ which satisfies (\hyperref[4A2]{4A2}) (resp. (\hyperref[5A2]{5A2}) and (\hyperref[5A1]{5A1}-(i))) such that $f$ must be unbounded from above. It can be seen that both Example \ref{ex4.6}(ii) and Example \ref{ex5.10}(ii) only need that $f$ is a positive constant.
\end{remark}

\begin{appendix}
\section*{Appendix}
\end{appendix}
\begin{appendix}
\section{Proofs of Lemma \ref{l2.3} and Lemma \ref{l2.4}}\label{appA}

\begin{proof}[Proof of Lemma \ref{l2.3}:] Since $g$ satisfies (\hyperref[2A1]{2A1}) and the RBSDE$(g,\xi,L)$ is dominated by $Y^1$, by \cite[Theorem 3.1]{EH}, the RBSDE$(g,\xi,L)$ admits a minimal solution $(y,z,k)$ (resp. a maximal solution $(Y,Z,K)$) such that $L_t\leq y_t\leq Y_t^1$ (resp. $L_t\leq Y_t\leq Y_t^1$). We now prove that $(y,z,k)$ is actually a minimal solution to the RBSDE$(g,\xi,L)$ such that $y\in\mathcal{C}_D.$ Let $(\tilde{Y},\tilde{Z},\tilde{K})$ be another solution to the RBSDE$(g,\xi,L)$ such that $y\in\mathcal{C}_D.$ In view of $y_t\wedge\tilde{Y}_t=\frac{1}{2}((y_t+\tilde{Y}_t)-|y_t-\tilde{Y}_t|)$, by applying Tanaka's formula to $y_t\wedge \tilde{Y}_t$, we have
\begin{eqnarray*}
&&y_t\wedge\tilde{Y}_t\\
&=&\xi+\frac{1}{2}\int_t^T((g(s,y_s,z_s)+g(s,\tilde{Y}_s,\tilde{Z}_s))-\textmd{sgn}(y_s-\tilde{Y}_s)(g(s,y_s,z_s)-g(s,\tilde{Y}_s,\tilde{Z}_s)))ds\\
&&+\frac{1}{2}((k_T+\tilde{K}_T)-(k_t+\tilde{K}_t))-\frac{1}{2}\int_t^T\textmd{sgn}(y_s-\tilde{Y}_s)d(\hat{K}_s-\tilde{K}_s)\\
&&-\frac{1}{2}\int_t^T((z_s+\tilde{Z}_s)-\textmd{sgn}(y_s-\tilde{Y}_s)(z_s-\tilde{Z}_s))\cdot dB_s
+\frac{1}{2}(\ell_T^0(y_s-\tilde{Y}_s)-\ell_t^0(y_s-\tilde{Y}_s)),
\end{eqnarray*}
where $\ell_t^0(y_s-\tilde{Y}_s)$ is the local time of $y_s-\tilde{Y}_s$ at time $t$ and level $0$. Observe that $dt\times dP$-\emph{a.e.,}
\begin{eqnarray*}
&&d(k_s+\tilde{K}_s)-\textmd{sgn}(y_s-\tilde{Y}_s)d(k_s-\tilde{K}_s)\\
&\geq&d(k_s+\tilde{K}_s)-|\textmd{sgn}(y_s-\tilde{Y}_s)dk_s
-\textmd{sgn}(y_s-\tilde{Y}_s)d\tilde{K}_s|\\
&\geq&d(k_s+\tilde{K}_s)-(|\textmd{sgn}(y_s-\tilde{Y}_s)|dk_s
+|\textmd{sgn}(y_s-\tilde{Y}_s)|d\tilde{K}_s)\\
&=&(1-|\textmd{sgn}(y_s-\tilde{Y}_s)|)d(k_s+\tilde{K}_s),
\end{eqnarray*}
it is not hard to check that the RBSDE$(g,\xi,L)$ is dominated by $y\wedge \tilde{Y}$. Then by \cite[Theorem 3.1]{EH} again, we get that the RBSDE$(g,\xi,L)$ admits a solution $(\hat{y},\hat{z},\hat{k})$ such that $L_t\leq \hat{y}_t\leq y_t\wedge\tilde{Y}_t\leq y_t\leq Y_t^1$, which implies that $y_t=\hat{y}_t,$ and thus $y_t\leq\tilde{Y}_t.$ The proof is complete.
\end{proof}

\begin{proof}[Proof of Lemma \ref{l2.4}:] Since $g$ satisfies (\hyperref[2A1]{2A1}) and the BSDE$(g,\xi)$ is dominated by $Y^2$ and $Y^1$, by \cite[Theorem 3.1]{EH}, we deduce that the BSDE$(g,\xi)$ admits a minimal solution $(y,z)$ (resp. a maximal solution $(Y,Z)$) such that $Y_t^2\leq y_t\leq Y_t^1$ (resp. $Y_t^2\leq Y_t\leq Y_t^1$). For another solution $(\tilde{Y},\tilde{Z})$ to the BSDE$(g,\xi)$ such that $\tilde{Y}_t\geq Y_t^2$, as discussed in the proof of Lemma \ref{l2.3}, we get that the BSDE$(g,\xi)$ is dominated by $Y^2$ and $\tilde{Y}_t\wedge y_t$, by \cite[Theorem 3.1]{EH} again, the BSDE$(g,\xi)$ admits a solution $(\hat{Y},\hat{Z})$ such that $Y_t^2\leq \hat{Y}_t\leq \tilde{Y}_t\wedge y_t\leq Y_t^1$. This implies that $(y,z)$ is a minimal solution to the BSDE$(g,\xi)$ such that $y_t\geq Y_t^2$.

We now prove that $(Y,Z)$ is a maximal solution to the BSDE$(g,\xi)$ such that $Y_t\leq Y^1_t$. Set $\tilde{g}(\cdot,\cdot,\cdot):=-g(\cdot,-\cdot,-\cdot)$. It can be checked that the BSDE$(\tilde{g},-\xi)$ is dominated by $-Y^1$ and $-Y^2$, and $(-{Y},-{Z})$ is a minimal solution to the BSDE$(\tilde{g},-\xi)$ such that $-Y^1_t\leq-Y_t\leq-Y^2_t$. Then by the argument above, we conclude that $(-{Y},-{Z})$ is a minimal solution to the BSDE$(\tilde{g},-\xi)$ such that $-Y_t\geq -Y^1_t$. For any solution $(\tilde{y},\tilde{z})$ to the BSDE$(g,\xi)$ such that $\tilde{y}_t\leq Y^1_t$, it is clear that $(-\tilde{y},-\tilde{z})$ is a solution to the BSDE$(\tilde{g},-\xi)$ such that $-\tilde{y}_t\geq -Y^1_t$, which implies $-\tilde{y}_t\geq -Y_t$. Thus $(Y,Z)$ is a maximal solution to the BSDE$(g,\xi)$ such that $Y_t\leq Y^1_t$.
\end{proof}
\section{A BMO property for quadratic BSDEs}\label{appB}
The following Lemma \ref{lmB.1} gives a BMO property for quadratic BSDEs, which slightly generalizes the corresponding results in \cite[Proposition 2.1]{BE} and \cite[Corollary 4.1]{B17}.
\begin{lemma}\label{lmB.1} Let the BSDE$(g,\xi)$ admit a solution $(Y,Z)$ such that $Y\in{\mathcal{S}}_D^{\infty}$. If there exist two nonnegative processes $\sqrt{\eta}\in {\mathcal{H}}^{BMO}_1$ and $C\in{\mathcal{S}}^\infty$, such that $g(t,Y_t,Z_t)\leq\eta_t+C_t|Z_t|^2$, $dt\times dP$-a.e., then $Z\in {\mathcal{H}}^{BMO}_d$.
\end{lemma}
\begin{proof} Set $M:=\sup_{t\in[0,T]}\|Y_t\|_\infty+\sup_{t\in[0,T]}\|C_t\|_\infty$ and for $n\geq 1$, $$\tau_n:=\inf\left\{t\geq0:\int_0^t|Z_s|^2ds>n\right\}\wedge T.$$
For $\beta>0$ and stopping time $\tau\leq\tau_n$, by applying It\^{o}'s formula to $e^{\beta Y_t}$, we have
\begin{align}
e^{\beta Y_\tau}+&\frac{\beta^2}{2}\int_\tau^{\tau_n}e^{\beta Y_s}|Z_s|^2ds\notag\\
&\leq e^{\beta Y_{\tau_n}}+\int_\tau^{\tau_n}\beta e^{\beta Y_s}\eta_sds+\int_\tau^{\tau_n}\beta e^{\beta Y_s}C_s|Z_s|^2ds-\int_\tau^{\tau_n}\beta e^{\beta Y_s}Z_s\cdot dB_s,\tag{b.1}\label{b.1}
\end{align}
which implies
$$\left(\frac{\beta^2}{2}-\beta M\right)E\left[\int_0^{\tau_n}e^{\beta Y_s}|Z_s|^2ds\right]
\leq e^{\beta M}+\beta e^{\beta M}E\left[\int_0^T\eta_sds\right].$$
Set $\beta>2M$. Since $\sqrt{\eta}\in {\mathcal{H}}^{BMO}_1$, by Fatou's Lemma, we have $Z\in{\mathcal{H}}^2_d.$ Then by setting $\tau_n=T$ in (\ref{b.1}), we conclude that, for each stopping time $\tau\leq T$,
$$\left(\frac{\beta^2}{2}-\beta M\right)E\left[\int_\tau^Te^{\beta Y_s}|Z_s|^2ds|{\mathcal{F}}_\tau\right]
\leq e^{\beta M}+\beta e^{\beta M}E\left[\int_\tau^T\eta_sds|{\mathcal{F}}_\tau\right].$$
Since $\sqrt{\eta}\in {\mathcal{H}}^{BMO}_1$, by setting $\beta>2M$, we have $Z\in {\mathcal{H}}^{BMO}_d$.
\end{proof}
\section{Some properties of the function $u_f(y)$}\label{appC}
Inspired by \cite{B17}, for $f(x)\in L_{loc}(D)$, $D=\mathbf{R}$ or $D=(0,\infty)$, we define
\begin{equation*}
u_f(y):=\int_1^y\exp\left(2\int^x_1f(z)dz\right)dx,\ \ y\in D.
\end{equation*}
The following properties for $u_f(x)$ come from \cite[Lemma A.1]{B17} and \cite[Lemma 2.1]{Zheng2}.
\begin{lemma}\label{lmC.1} For $f(x)\in L_{loc}(D)$, the following properties of $u_f(x)$ hold:

(i) $u_f(x)\in W^2_{1,loc}(D)$, in particular, $u_f(x)\in C^1(D)$;

(ii) $u_f(x)$ is strictly increasing, and when $f(\cdot)\geq0$ (resp. $\leq0$), $u_f(x)$ is convex (resp. concave);

(iii) $u_f''(x)-2f(x)u_f'(x)=0,\ a.e.$ on $D$;

(iv) $u_f^{-1}(x)\in W^2_{1,loc}(V)$, where $V:=\{y:y=f(x),x\in D\}$. Moreover, $u_f^{-1}(x)\in C^1(V)$ and is strictly increasing;

(v) If $l\in L_{1,loc}(D)$ and $l(x)\leq f(x), a.e.$, then for each $x\in D,$ $u_l(x)\leq u_f(x)$.
\end{lemma}

\begin{example}\label{exC.2} We show some typical examples of $u_f(x)$:

(i) If $f(x)=\beta, x\in\mathbf{R}$ for some constant $\beta$, then
\begin{equation*}
u_\beta(x):=u_f(x)=\left\{
      \begin{array}{ll}
        x-1, & \beta=0; \\
        \frac{1}{2\beta}(\exp(2\beta(x-1))-1), & \beta\neq0,
      \end{array}
    \right.\quad x\in\mathbf{R};
\end{equation*}

(ii) If $f(x)=\frac{\beta}{x}, x\in(0,\infty)$ for some constant $\beta$, then
\begin{equation*}
u_{\frac{\beta}{\cdot}}(x):=u_f(x)=\left\{
      \begin{array}{ll}
        \frac{1}{1+2\beta}(x^{1+2\beta}-1), & \beta\neq-\frac{1}{2}; \\
        \ln(x), & \beta=-\frac{1}{2},
      \end{array}
    \right.\quad x\in(0,\infty).
\end{equation*}
\end{example}
\section{On the convexity of the function $y^{-r}|z|^2$}\label{appD}

\begin{lemma}\label{lmD.1} $y^{-r}|z|^2$ is convex on $(0,\infty)\times\mathbf{R}^d$ if and only if $0\leq r\leq1.$
\end{lemma}

\begin{proof}
For $y\in(0,\infty)$ and $z=(z_1,\cdots,z_d)^\textsf{T}\in\mathbf{R}^d$, the Hessian matrix of $y^{-r}|z|^2$ is:
\begin{equation*}
H(y,z)=\begin{bmatrix}
r(r+1)y^{-r-2}|z|^2 & -2r y^{-r-1} z_1 & \cdots & -2r y^{-r-1} z_d \\
-2r y^{-r-1}z_1 & 2y^{-r} & \cdots & 0 \\
\vdots & \vdots & \ddots & \vdots \\
-2r y^{-r-1}z_d & 0 & \cdots & 2y^{-r}
\end{bmatrix}
\end{equation*}
It is clear that $y^{-r}|z|^2$ is convex on $(0,\infty)\times\mathbf{R}^d$ if and only if $H(y,z)$ is positive semi-definite on $(0,\infty)\times\mathbf{R}^d$, i.e., for all $v\in\mathbf{R}^{d+1}$, we have $v^\textsf{T}H(y,z)v\geq0$ on $(0,\infty)\times\mathbf{R}^d.$

For all $v =(v_1,v_2)^\textsf{T}=(v_1, v_{2,1},\cdots,v_{2,d})^\textsf{T}\in\mathbf{R}^{d+1}$, we have
\begin{equation*}\label{d.1}
v^\textsf{T}H(y,z)v=y^{-r}(r(r+1)|z|^2y^{-2}v_1^2-4ry^{-1}v_1(v_2\cdot z)+2|v_2|^2).\tag{d.1}
\end{equation*}
Clearly, if $|z|=0,$ then for all $v$, we have $v^\textsf{T}H(y,0)v\geq0.$ Thus, we only need to consider the case that $|z|\neq0.$ In this case, by setting $a:=y^{-1}|z|$, $\hat{z}:=|z|^{-1}z$, $b:=v_2\cdot \hat{z}$, and $w:=v_2-b\hat{z}$, we get from (\ref{d.1}) that
\begin{align*}
v^\textsf{T}H(y,z)v&=y^{-r}(r(r+1)a^2v_1^2-4rav_1b+2b^2+2|w|^2)\\
&\geq y^{-r}(r(r+1)a^2v_1^2-4rav_1b+2b^2)\\
&=y^{-r}(2(b-rav_1)^2+r(1-r)a^2v_1^2)
\end{align*}
From this, we get that if $0\leq r\leq1,$ then $y^{-r}|z|^2$ is convex on $(0,\infty)\times\mathbf{R}^d$. Conversely, if $y^{-r}|z|^2$ is convex on $(0,\infty)\times\mathbf{R}^d$, by setting $v_1=1$ and $v_2=ry^{-1}z$, then we get from (\ref{d.1}) that $v^\textsf{T}H(y,z)v=y^{-2-r}r(1-r)|z|^2$, which implies $0\leq r\leq1.$ The proof is complete.
\end{proof}
\section{An existence and uniqueness result for a special RBSDE}\label{appE}
\begin{lemma}\label{lmE.1} Let $D=(0,\infty)$, $\gamma\in{\mathcal{C}}$, and $l\in{\mathcal{C}}_D$. Let $\kappa$ and $\nu$ be two nonnegative constants. If $u_{\frac{-\nu}{\cdot}}(e^{-\int_0^T\gamma_sds}\xi)\in L^p(\mathcal{F}_T)$ and $u_{\frac{-\nu}{\cdot}}(e^{-\int_0^\cdot\gamma_sds}l)\in \mathcal{S}^p$ such that $\xi\leq l_T$, then the RBSDE$(-\gamma_ty+\kappa|z|-\frac{\nu}{y}|z|^2,-\xi,-l_t)$ admits a unique solution $(Y,Z,K)$ such that $-Y\in{\mathcal{C}}_D$ and $u_{\frac{-\nu}{\cdot}}(-e^{-\int_0^\cdot\gamma_sds}Y)\in\mathcal{S}^p$.
\end{lemma}

\begin{proof} By the definition of $u_{\frac{-\nu}{\cdot}}$ (see Appendix \ref{appC}), we have
\begin{equation*}\label{e.1}
-u_{\frac{-\nu}{\cdot}}(-y)=\int_{-1}^y\exp\left(2\int_{-1}^x-\frac{\nu}{z}dz\right)dx,\ \ y\in(-\infty, 0).\tag{e.1}
\end{equation*}
It follows from \cite[Corollary 3.5]{Zheng3} and (\ref{e.1}) that the RBSDE$(\kappa|z|-\frac{\nu}{y}|z|^2,-e^{-\int_0^T\gamma_sds}\xi,-e^{-\int_0^t\gamma_sds}l_t)$ admits a unique solution $(\hat{Y},\hat{Z},\hat{K})$ such that $-\hat{Y}\in{\mathcal{C}}_D$ and $u_{\frac{-\nu}{\cdot}}(-\hat{Y})\in \mathcal{S}^p$. Set
\begin{equation*}
Y_t:=e^{\int_0^t\gamma_sds}\hat{Y}_t,\ Z_t:=e^{\int_0^t\gamma_sds}\hat{Z}_t,\ K_t:=\int_0^te^{\int_0^s\gamma_rdr}d\hat{K}_s,\ \ t\in[0,T].
\end{equation*}
By It\^{o}'s formula, we can deduce that $(Y,Z,K)$ is a unique solution to the RBSDE$(-\gamma_ty+\kappa|z|-\frac{\nu}{y}|z|^2,-\xi,-l)$ such that $-Y\in{\mathcal{C}}_D$ and $u_{\frac{-\nu}{\cdot}}(-e^{-\int_0^\cdot\gamma_sds}Y)\in\mathcal{S}^p$.
\end{proof}
\end{appendix}

\ \\
\textbf{Declaration of competing interest}\\

The authors declare that they have no known competing financial interests or personal
relationships that could have appeared to influence the work reported in this paper.

\ \\
\textbf{Acknowledgement}\\

The author is very grateful to the anonymous referee for many valuable suggestions. The author would like to thank Prof. Shengjun Fan for helpful discussions.

\end{document}